\documentclass[11pt]{article}
\usepackage{amssymb,amsmath,amsfonts,amsthm,mathtools}
\usepackage{bbm}
\usepackage{bm} 
\usepackage[margin=1in]{geometry}
\usepackage{setspace}
\usepackage{enumitem}
\usepackage{xcolor,graphicx,subcaption}
\usepackage{fancyhdr,lastpage}
\usepackage[utf8]{inputenc}
\usepackage[numbers,sort&compress]{natbib}
\usepackage[unicode=true,pdfusetitle,bookmarks=false, breaklinks=true,pdfborder={0 0 0},pdfborderstyle={},backref=false,colorlinks=false]
{hyperref}
\usepackage{url}

\usepackage{amsthm}
\newtheorem{theorem}{Theorem}[section]
\newtheorem{lemma}[theorem]{Lemma}

\newtheorem{definition}[theorem]{Definition}
\newtheorem*{assumptions*}{Assumptions}

\newtheorem{manualtheoreminner}{Lemma}
\newenvironment{manualtheorem}[1]{%
  \renewcommand\themanualtheoreminner{#1}%
  \manualtheoreminner
}{\endmanualtheoreminner}

\DeclareMathAlphabet\mathbfcal{OMS}{cmsy}{b}{n}
\newcommand{\balpha}{\bm{\alpha}}
\newcommand{\bbeta}{\bm{\beta}}

\newcommand{\bpsi}{\bm{\psi}}
\newcommand{\e}{\varepsilon}
\newcommand{\x}{\mathbf{x}}
\newcommand{\R}{{\mathbb{R}}}

\newcommand{\A}{{\mathcal{A}^{(t_0,t_f)}}}
\newcommand{\ASig}{{\mathcal{A}_{\Sigma}^{(t_0,t_m)}}}
\newcommand{\ASigO}{{\overline{\mathcal{A}}_\Sigma^{(t_0,t_m)}}}
\newcommand{\I}{\mathcal{I}}

\newcommand{\Isig}{{I^{(t_0,t_f)}}}

\newcommand{\IsigMol}{{I_{\e}^{(t_0,t_f)}}}
\newcommand{\Ilim}{{\overline{I}^{(t_0,t_f)}}}

\newcommand{\ISigO}{{\overline{I}_{\Sigma}^{(t_0,t_m)}}}

\newcommand{\fscaled}{{\mathbf{F}^{\e}}}
\newcommand{\wlim}{\,\overset{H^1}{\rightharpoonup}\,}
\DeclareMathOperator*{\argmin}{arg\,min}
\newcommand{\B}{\overline{B}}

\newcommand*\samethanks[1][\value{footnote}]{\footnotemark[#1]}

\graphicspath{{FiguresSmaller/}}

\title{Most probable transition paths in piecewise-smooth stochastic differential equations}
\author{Kaitlin Hill\thanks{Department of Mathematics, St.~Mary's University, San Antonio, TX 78228, USA}~\,\thanks{Corresponding author, \textit{Email}: khill5@stmarytx.edu}\,, Jessica Zanetell\thanks{Department of Mathematics, Wake Forest University, Winston-Salem, NC 27109, USA}\,, and John A.~Gemmer\samethanks[3]~\,\thanks{\textit{Email}: gemmerj@wfu.edu}}
\date{}

\begin{document}
\maketitle

\begin{abstract}
\noindent We develop a path integral framework for determining most probable paths for a class of systems of stochastic differential equations with piecewise-smooth drift and additive noise. This approach extends the Freidlin-Wentzell theory of large deviations to cases where the system is piecewise-smooth and may be non-autonomous. In particular, we consider an $n-$dimensional system with a switching manifold in the drift that forms an $(n-1)-$dimensional hyperplane and investigate noise-induced transitions between metastable states on either side of the switching manifold. To do this, we mollify the drift and use $\Gamma-$convergence to derive an appropriate rate functional for the system in the piecewise-smooth limit. The resulting functional consists of the standard Freidlin-Wentzell rate functional, with an additional contribution due to times when the most probable path slides in a crossing region of the switching manifold. We explore implications of the derived functional through two case studies, which exhibit notable phenomena such as non-unique most probable paths and noise-induced sliding in a crossing region.
\end{abstract}

\noindent{\footnotesize \textbf{Keywords:} Piecewise smooth dynamical systems, Filippov systems,  Freidlin-Wentzell rate functional, Gamma-convergence, noise induced tipping, rare events

\noindent \textbf{2000 MSC:} 37H10, 37J45}


\section{Introduction}
\label{sec:intro}
In dynamical systems, a \emph{tipping event} is loosely defined as occurring when a sudden or small change to a variable or parameter induces a large change to the state of the system. While tipping is often studied within the context of climate applications \cite{berglund2002metastability, Eisenman09, Wieczorek10, eisenman2012factors, rothman2019characteristic, arnscheidt2020routes, lohmann2021risk}, it also has broad applications in  ecology \cite{vanselow2019very, o2020tipping}, ecosystems \cite{drake2010early, scheffer2008pulse} and epidemiology \cite{forgoston2011maximal, schwartz2011converging, Billings2018, hindes2018rare}, to name a few. While there is no precise mathematical definition of a tipping event, in \cite{ashwin2012tipping} it was proposed that tipping events could be classified according to whether the underlying mathematical mechanism involves, predominantly, a bifurcation (B-tipping), noise-induced transitions (N-tipping), or transitions between basins of attraction induced by fast changes in parameters (R-tipping), i.e. rate-induced tipping; see also \cite{thompson2011climate,lenton2011early, halekotte2020minimal, alkhayuon2021phase}. 

The term `noise-induced tipping' encompasses a range of phenomena such as transitions between metastable states, bursting,  stochastic resonance, stochastic coherence and stochastic synchronization \cite{lindner2004effects} and can be analyzed from several points of view including the Freidlin-Wentzell (FW) theory of large deviations, \cite{freidlin2012random, ren2004minimum, forgoston2018primer}, the path integral framework \cite{durr1978onsager, chaichian2001path}, transition path theory \cite{vanden2010transition, finkel2020path}, and formal asymptotics \cite{maier1997limiting, muratov2008noise}. In particular, for smooth autonomous dynamical systems additively perturbed by Gaussian white noise, the FW theory provides a framework for computing most probable transitions as minimizers of a rate functional in the asymptotic limit of vanishing noise strength. The benefit of the FW framework is that most probable transition paths can be numerically computed using iterative schemes such as the string method  for gradient systems \cite{weinan2002string}, the minimum action method \cite{ren2004minimum}, the geometric minimum action method \cite{heymann2008geometric} and explicit gradient descent \cite{lindley2013iterative}. Moreover, knowledge of the most probable transition path can then be coupled with statistical techniques such as importance sampling to compute quantities of interest such as the expected time of tipping for nonzero noise \cite{donovan2011iterative, forgoston2018primer}. 

In this paper we extend the path integral framework to a class of differential equations with additive noise and piecewise-smooth drift. There has been considerable recent interest in understanding the dynamics of piecewise-smooth stochastic dynamical systems, particularly in relation to rare events and tipping, or transitions between metastable states, but also in varying applications \cite{simpson2018influence,staunton2020discontinuity,staunton2020estimating}. In \cite{chen2013weak,chen2014large}, Chen et al.~use the backward Fokker-Planck technique to derive the distribution of paths in a model of dry friction. Notably, in \cite{chen2013weak}, Chen et al.~use a similar method to that of the present study, smoothing out the piecewise-smooth vector field to numerically analyze the desingularized SDE. In \cite{baule2010stick,baule2009path}, Baule et al.~use the path integral approach, as in the present study, to investigate most probable paths in a piecewise-smooth model of stick-slip friction. Beyond mechanical systems, the study of noise-induced tipping in piecewise-smooth systems has applications in biology \cite{simpson2011mixed,piltz2014prey} and climate models \cite{serdukova2017metastability, moon2017stochastic, yang2020tipping, morupisi2021analysis,monahan2002stabilization}. To our knowledge, most probable paths in general stochastic piecewise-smooth systems have not yet been addressed.

\subsection{Most probable paths in smooth systems}

During rare event transitions from one metastable state to another, a stochastic system will follow paths according to some distribution, which in the limit of vanishing noise strength is generally singly-peaked along a most probable transition path. Following \cite{freidlin2012random}, for a smooth vector field $\mathbfcal{F}$, the most probable transition path between $\x_0$ and $\x_f$ can be defined as follows. We first define an admissible set of transition paths $\A$ by
\[
\A=\{\balpha \in H^{1}([t_0,t_f];\R^n):\balpha(t_0)=\x_0 \text{ and } \balpha(t_f)=\x_f\},
\]
where $H^1([t_0,t_f];\R^n)$ is the Hilbert space of weakly differentiable curves $\balpha$ satisfying $\balpha\in H^1$ if and only if $ \int_{t_0}^{t_f}|\balpha|^2+|\dot{\balpha}|^2\,dt<\infty$. The most probable transition path $\balpha^*\in \A$ is then defined to be the global minimizer of the Freidlin-Wentzell rate functional $\Isig:\A\mapsto \overline{\R},$ given by
\begin{equation}\label{eq:FW-functional}
\Isig[\balpha]=\int_{t_0}^{t_f}\left\| \dot{\balpha}(t)-\mathbfcal{F}(\balpha(t),t)\right\|^2\,dt,
\end{equation}
so that $\balpha^*=\argmin_{\balpha\in \A} \Isig[\balpha]$. For this functional, the necessary condition satisfied by second differentiable  minimizers is the following Euler-Lagrange equations \cite{forgoston2018primer}: 
\[
\begin{aligned}
\ddot{\balpha} = \mathbfcal{F}_t + \dot{\balpha}\left( \nabla\mathbfcal{F} - \nabla\mathbfcal{F}^{\intercal} \right) + \nabla\mathbfcal{F}^{\intercal}\mathbfcal{F}.
\end{aligned}
\]
Here, the subscript $t$ refers to the partial derivative with respect to time. 

\subsection{Framing the piecewise-smooth problem}

We determine the most probable transition path of a noise-induced tipping event when the drift is piecewise-smooth. In general, for systems with a piecewise-smooth drift, the appropriate rate functional to minimize is not known. Minimizers of the Freidlin-Wentzell rate functional may not be well-defined when a region of the vector field is not continuous, or even Lipschitz continuous. It is not clear how the most probable path may traverse across a switching manifold, which may have an attracting or repelling sliding region that introduces more complex dynamics. Specifically, letting $\mathbf{F}^{\pm}:\R^n\mapsto \R^n$ be smooth vector fields, we consider a system of stochastic differential equations of the form
\begin{equation}\label{eq:Intro:SDE}
d\x_t=\mathbf{F}(\x_t)dt+\sigma d\mathbf{W}_t,
\end{equation}
where $\x=(x,\mathbf{y})\in \R^n$ such that $x\in\R$ and $\mathbf{y}\in\R^{n-1}$, $\sigma\in \R$, $\mathbf{W}=(W_1,\ldots W_n)$ is an $n$-dimensional Wiener process, and $\mathbf{F}:\R^n\setminus \{x=0\}\mapsto \R^n$ is defined by
\begin{equation}\label{eq:Intro:PWVF}
\mathbf{F}(\x)=\begin{cases}
\mathbf{F}^{+}(\x), &x>0,\\ 
\mathbf{F}^{-}(\x), & x<0.
\end{cases}
\end{equation}
We let $S_+=\{\x\in \R^n: x>0\}$, $S_{-}=\{\x\in \R^n: x<0\}$, and $\Sigma=\{\x\in \R^n: x=0\}$. The set $\Sigma$ is often called the \textit{switching manifold} or \textit{discontinuity boundary} \cite{diBernardo08}. In general, $\Sigma$ may be defined as the zero level set of a smooth function $H:\R^n\mapsto\R,$ but we assume for simplicity that $\Sigma$ is a hyperplane in $\R^{n-1}$.

We define the deterministic skeleton of Equation \eqref{eq:Intro:SDE} as the dynamical system
\begin{equation}\label{eq:Intro:DS}
\dot{\x}=\mathbf{F}(\x). 
\end{equation}
We assume that there exist asymptotically stable fixed points $\x_0\in S_+$ and $\x_f \in S_-$ of Equation \eqref{eq:Intro:DS}, such that $\mathbf{F}^+(\x_0)=0$ and $\mathbf{F}^-(\x_f)=0$. By noise-induced transitions, we mean realizations of Equation \eqref{eq:Intro:SDE} that transition from the basin of attraction of $\x_0$ to that of $\x_f$. More precisely, we define a \textit{noise-induced transition} on the interval $[t_0,t_f]$ as a solution to the stochastic boundary value problem given by realizations of Equation \eqref{eq:Intro:SDE} that satisfy $\x(t_0)=\x_0$ and $\x(t_f)=\x_f$. Here the boundary conditions mean that the process is conditioned to transition from $\x_0$ to $\x_f$.

\subsection{Background on piecewise-smooth systems}
\label{sec:PWS-background}

Dynamics that occur entirely within the smooth regions $S_{\pm}$ can be fully described by regular (smooth) dynamical systems theory. On the other hand, dynamics on a switching manifold $\Sigma$ may not be defined and are typically imposed \textit{a posteriori}, e.g.~using Filippov's convex combination \cite{diBernardo08,diBernardo,filippov2013differential}. In general, the lack of smoothness requirements across the switching manifold allows for more diverse phenomena than in smooth systems of similar dimension, since the vector field across $\Sigma$ may be discontinuous or even point in opposite directions. Regions of $\Sigma$ where either occurs are called \textit{sliding regions}. A solution that reaches $\Sigma$ at a sliding region may ``slide'' along it. On the other hand, \textit{crossing regions} occur on $\Sigma$ where sliding is not possible. We differentiate these regions using the following notation:
\begin{enumerate}
	\item \emph{Attracting sliding regions}: $\Sigma_{A}=\{\x\in \Sigma: F_1^+(0,\mathbf{y})\leq 0,\,F_1^-(0,\mathbf{y})\geq0\}$,
	\item \emph{Repelling sliding regions}: $\Sigma_{R}=\{\x\in \Sigma: F_1^+(0,\mathbf{y})\geq0,\, F_1^-(0,\mathbf{y})\leq 0\},$
	\item \emph{Positive crossing regions}: $\Sigma_{+}=\{\x\in \Sigma: F^{\pm}_1(0,\mathbf{y})>0\}$,
	\item \emph{Negative crossing regions}: $\Sigma_{-}=\{\x\in \Sigma: F^{\pm}_1(0,\mathbf{y})<0\} $,
\end{enumerate}
where ${F}^{\pm}_{1}$ is the first component of $\mathbf{F}^{\pm}$.

When solutions may slide along the switching manifold $\Sigma,$ we impose a flow using the Filippov convex method \cite{filippov2013differential, diBernardo08}. That is, we define the \textit{sliding flow} as a convex combination of $\mathbf{F}^{\pm},$
\begin{equation}\label{eq:sliding-flow}
\mathbf{F}^{s}(0,\mathbf{y}) = \lambda \mathbf{F}^+(0,\mathbf{y}) + (1-\lambda)\mathbf{F}^-(0,\mathbf{y})
\end{equation}
with 
\begin{equation}\label{eq:gamma}
\lambda(\mathbf{y}) = \frac{F_1^-(0,\mathbf{y})}{F_1^-(0,\mathbf{y}) - F_1^+(0,\mathbf{y})}\in[0,1].
\end{equation}
Notice that Equation \eqref{eq:gamma} fixes $\lambda(\mathbf{y})$. 
This flow naturally arises in the context of our main result without \textit{a priori} imposing it; see Theorem \ref{thm:gamma-convergence}. Also, note that $\mathbf{F}^{s}(0,\mathbf{y})$ may vanish for some point $(0,\mathbf{y})\in\Sigma$. Since it is not an equilibrium of the smooth flow, we say that such a point $\x_s=(0,\mathbf{y}_s)\in\Sigma$ is a \textit{pseudoequilibrium} if it is an equilibrium of the sliding flow; i.e., for some $\lambda\in (0,1),$ $\mathbf{F}^{s}(0,\mathbf{y}_s)=\mathbf{0}$.

Filippov's convex combination interprets solutions of System \eqref{eq:Intro:PWVF},\eqref{eq:Intro:DS} as a continuous curve $\x(t)$ satisfying the following conditions:
\begin{equation}\label{eq:deterministic:odeFillipov}
\begin{cases}
\dot{\x}(t)=\mathbf{F}(\x),\\
\x(0)=\x_0,
\end{cases}
\end{equation}
where $\mathbf{F}:\R^n \mapsto \R^n$ is defined as
\begin{equation}\label{eq:deterministic:odeFillipovField}
\mathbf{F}(\x)=\begin{cases}
\mathbf{F}^+(\x), & \x\in S_+ \cup \Sigma_{+},\\
\mathbf{F}^-(\x), & \x \in S_{-} \cup  \Sigma_{-},\\
\mathbf{F}^s(\x), & \x \in \Sigma_{A}\cup \Sigma_{R}.
\end{cases}
\end{equation} 
That is, for points in time in which a solution curve $\x(t)$ intersects $\Sigma$ it will either cross $\Sigma,$ tracking the flow of either $\mathbf{F}^+$ or $\mathbf{F}^-$ depending on whether the curve is crossing into $S_+$ or $S_-$, or it will track $\Sigma$ in the sliding regions until entering a crossing region. Note, implicit in this definition is that solution curves are differentiable everywhere except on $\Sigma$. 

\begin{figure}[t!]
	\centering
	\begin{subfigure}{.2\textwidth}
		\centering
		\includegraphics[scale=1]{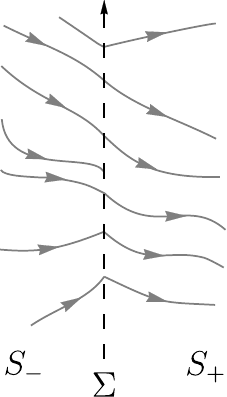}
		\caption{}
	\end{subfigure}%
	\begin{subfigure}{.2\textwidth}
		\centering
		\includegraphics[scale=1]{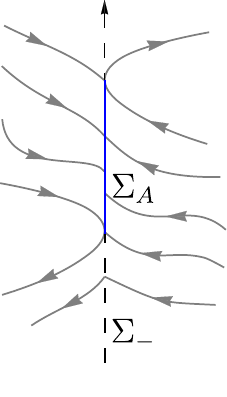}
		\caption{}
	\end{subfigure}%
	\begin{subfigure}{.2\textwidth}
		\centering
		\includegraphics[scale=1]{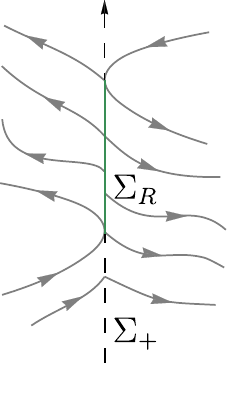}
		\caption{}
	\end{subfigure}%
	\begin{subfigure}{.35\textwidth}
		\centering
		\includegraphics[scale=1]{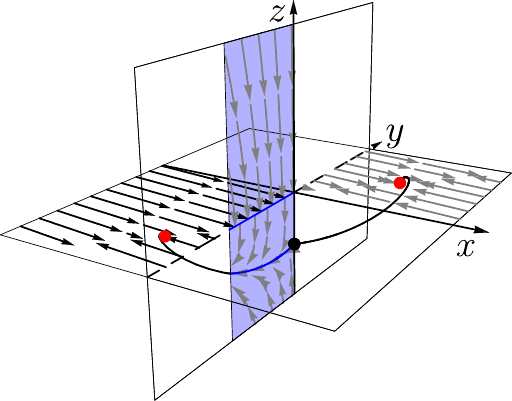}
		\caption{}
	\end{subfigure}
	\caption{(a)-(c) Phase plane of a generic $\mathbf{F}$ with $n=2$, with (a) a crossing region, (b) crossing and attracting sliding, and (c) crossing and repelling sliding. (d) Phase diagram of a piecewise-smooth version of the Lorenz 63 model, System \eqref{eq:LorenzFilippov}. Red points are equilibria of $S_{\pm},$ and the black point is a pseudoequilibrium of the imposed flow on $\Sigma$. The sliding region, $\Sigma_A,$ is shaded blue. The black curve connects the two equilibria, going through the pseudoequilibrium at the origin and sliding in $\Sigma_A$. The blue curve is the sliding trajectory. Parameters used are $\sigma_+ = 10$, $\beta_+ = 2$, $\rho_{\pm} = 2$, $\sigma_- = 11$, $\beta_- = 3$, $x_0 = 1$, $y_0 = -1$, and $z_0 = 0$.}
	\label{fig:generic-switching}
\end{figure}

To provide a brief illustration, Figure \ref{fig:generic-switching}(a)-(c) shows some possible dynamics across $\Sigma$ in a two-dimensional system, with (a) $\Sigma_+$, (b) $\Sigma_{\pm}$ and $\Sigma_A$, and (c) $\Sigma_{\pm}$ and $\Sigma_R$. In (a) the flow traverses across $\Sigma$ continuously but not smoothly; in (b) and (c) the flow may also remain on $\Sigma,$ in intervals indicated by the blue line (for $\Sigma_A$) and the green line (for $\Sigma_R$). Note that the flow leaving a repelling sliding region is non-unique in forward time, while the flow leaving an attracting sliding region is non-unique in backward time. In addition, complex dynamics may appear on $\Sigma$ itself, depending on the manner in which the imposed vector field on $\Sigma$ transitions between $S_+$ and $S_-$. Figure \ref{fig:generic-switching}(d) shows one example of a path traversing from a stable equilibrium in $S_+$, through the switching manifold in an attracting sliding region using dynamics imposed using Filippov's convex combination \cite{filippov2013differential}, then traversing to a stable equilibrium in $S_-$. The system in this illustration is a piecewise-smooth version of the Lorenz 63 model,
\begin{equation}\label{eq:LorenzFilippov}
\begin{pmatrix} \dot{x} \\ \dot{y} \\ \dot{z} \end{pmatrix} =     
\begin{pmatrix}
\sigma_-\left(y-y_0-x+x_0\right) \\
(x-x_0)(\rho_--z+z_0)-y+y_0 \\
(x-x_0)(y-y_0) - \beta_- (z-z_0)
\end{pmatrix},\; x<0,\mbox{ and}\;
\begin{pmatrix} \dot{x} \\ \dot{y} \\ \dot{z} \end{pmatrix} =
\begin{pmatrix}
\sigma_+\left(y-x\right) \\
x(\rho_+-z)-y \\
xy - \beta_+ z
\end{pmatrix},\; x>0,
\end{equation}
originally studied in \cite{Lorenz63}. Here, $\sigma_{\pm},\rho_{\pm},\beta_{\pm}>0$ are parameters and $x_0,y_0,z_0\in\R$. For this example, $S_{\pm}$ and $\Sigma$ are defined as previously, with $n=3$.

\subsection{Outline of analysis}

For our problem, noise-induced transitions must cross $\Sigma$ and the Freidlin-Wentzell rate functional, Equation \eqref{eq:FW-functional}, must be modified to account for the discontinuity in $\mathbf{F}$ and possible discontinuity-induced dynamics. 
For example, a most probable path that reaches a repelling sliding region $\Sigma_R$ leaves it non-uniquely due to the non-uniqueness of the flow, as illustrated in the case studies in Sections \ref{sec:case-study-2D} and \ref{sec:case-study-non-autonomous}. Additionally, if one were to naively piece together the most probable paths of the smooth regions and neglect possible contribution from the manner in which the path crosses the switching manifold, this may lead to paths that do not reflect a global minimum of the rate functional.

We resolve these issues by smoothing out $\mathbf{F}$ via mollification with a compactly supported, smooth, radially symmetric kernel of characteristic width $\e$, and consider the sequence of minimizers to the rate functional Equation \eqref{eq:FW-functional} as $\e\rightarrow 0$. We show for any piecewise-smooth system of the form given by System \eqref{eq:Intro:SDE},\eqref{eq:Intro:PWVF} that the most probable path minimizes Equation \eqref{eq:FW-functional} in the smooth regions $S_{\pm}$, plus an additional functional whose contribution represents time spent sliding in a crossing region of the switching manifold. For $\mathbf{F} = (F_1, \mathbf{G})$ such that $F_1:\R^n\setminus \Sigma\mapsto \R$ and $\mathbf{G}:\R^n\setminus \Sigma\mapsto \R^{n-1},$ and $\balpha=(\alpha,\bbeta)\in\A$ and $0<\e\ll 1$, the appropriate rate functional with both contributions is
\[
\begin{aligned}
\Ilim[\balpha]=&\int_{\{t\in S_+\cup S_-\}}\|\dot{\balpha}(t)-\mathbf{F}(\balpha(t))\|^2\,dt \\
&+\int_{\{t\in\Sigma_{\pm}\}}\min_{\lambda\in [0,1]}\left\{
\left[ \lambda F_1^+(0,\bbeta) + (1-\lambda)F_1^-(0,\bbeta)\right]^2 \right. \\
&\left. \mathmakebox[\widthof{$ +\int_{\I_{\Sigma}[\balpha]}\min\qquad$ }][l]{}  + \left| \dot{\bbeta} - \lambda\mathbf{G}^+(0,\bbeta) - (1-\lambda)\mathbf{G}^-(0,\bbeta)\right|^2
\right\}\,dt.
\end{aligned}
\]
We make the following observations about this result:
\begin{enumerate}
	\item The form of the above rate functional is independent of the chosen mollifier.
	\item When there is a sliding region, the most probable path may track the Filippov dynamics; that is, no additional contribution is made during times where the most probable path slides via the imposed flow on the switching manifold.
\end{enumerate}

We rigorously prove the above result using a technique from the calculus of variations called $\Gamma$-convergence to show that $\Ilim$ is the $\Gamma$-limit of the sequence of rate functionals for the mollified system \cite{dal2012introduction, braides2002gamma}. The $\Gamma$-limit is the natural notion of a limiting functional in this context in the sense that minimizers of the rate functional with the mollified drift converge (weakly) to a minimizer of $\Ilim$ in the limit as $\varepsilon\rightarrow 0$. While $\Gamma$-convergence has been used to study the convergence of minimizers of the Onsager-Machlup functional in the small noise limit \cite{pinski2012gamma, lu2017gaussian, li2021gamma,  ayanbayev2021gammaI, ayanbayev2021gammaII}, to our knowledge we are the first to use $\Gamma$-convergence to compute an appropriate extension of the Freidlin-Wentzell rate functional for piecewise smooth systems. 

The utility of the derived rate functional is demonstrated through two case studies. Each case study also presents distinct phenomena in their most probable transition paths that are not possible in systems with smooth drift. The first case study analyzes the most probable transition path in a two-dimensional piecewise-linear system, a simple demonstrative case in which the typical rate functional formulation of the most probable path breaks down. Notably, in this case the most probable path follows the switching manifold along an attracting sliding region, but it does not follow the switching manifold along a repelling sliding region.

The second case study analyzes the most probable path in a simple one-dimensional periodically-forced piecewise-smooth system, constructed similarly to some conceptual models for Arctic energy balance \cite{eisenman2012factors,hill2016analysis}. A significant portion of the analysis of the deterministic dynamics and the Monte Carlo simulations were first derived in Zanetell's Master's thesis \cite{Zanetell}. Here we again observe the breakdown of the typical formulation of the most probable path. Additionally, we observe an emergent dynamical phenomenon, which we describe as noise-induced sliding, in which the most probable path follows the switching manifold in a crossing region, despite the energy cost of the associated rate functional.

In both case studies, we note in some regimes the most probable path(s) predicted by our result do not match observed Monte Carlo simulations. Several factors may contribute to this discrepancy, including needing a smaller noise strength $\sigma,$ a higher-order functional in the smoothing parameter $\e,$ or including the correction term in the Onsager-Machlup functional. Such possibilities will be explored in future work.

The paper is organized as follows: in Section \ref{sec:notation} we briefly highlight notation we use and some technical assumptions on $\mathbf{F}$. Then in Section \ref{sec:most-probable-paths} we derive the mollified vector field and use $\Gamma-$convergence to determine the appropriate limiting functional for calculating the most probable path for stochastic differential equations of the form given by System \eqref{eq:Intro:SDE},\eqref{eq:Intro:PWVF}. Next, in Section \ref{sec:case-study-2D} we provide a simple case study of a planar linear piecewise-smooth system as an illustration of the most probable path calculated in Section \ref{sec:most-probable-paths}. Finally, in Section \ref{sec:case-study-non-autonomous} we provide a case study that extends the most probable path derivation to the case of a one-dimensional non-autonomous system.

\section{Notation and Assumptions}
\label{sec:notation}

In this paper we use the following conventions for notation:
\begin{enumerate}
	\item We use the convention that $\| \cdot \|$ represents the $L^2$ norm for functions and $|\cdot |$ represents the $\ell^2$ norm for vectors in $\R^n$. Furthermore, we use $\langle \cdot , \cdot \rangle$ to denote the inner product in $\ell^2$.
	
	\item Subscripts refer to the index of a vector unless otherwise indicated. E.g., for $\mathbf{F}=(F_1,\,F_2,\,\ldots,\,F_n),$ we write the components as $F_j,$ $j=1,\,\ldots,\,n$.
	
	\item Unless noted otherwise, $\mathbf{F}$ may be non-autonomous. Although our analysis applies to non-autonomous systems, for simplicity of presentation we suppress the time dependence.
	
	\item Since we assume that the underlying drift is only piecewise-smooth in the first component, we use a convenient shorthand that separates this first component from the remaining ones. Namely, we write $\x=(x,\mathbf{y})\in\R^n,$ where $x\in\R$ and $\mathbf{y}\in\R^{n-1},$ and $\mathbf{F} = (F_1, \mathbf{G}),$ where $F_1:\R^n\setminus \{x=0\}\mapsto \R$ and $\mathbf{G}:\R^n\setminus \{x=0\}\mapsto \R^{n-1}$. 
	
	\item For the path $\balpha\in\A,$ we separate the first component from the remaining ones as $\balpha=\left(\alpha,\bbeta\right),$ so that $\alpha\in H^1\left([t_0,t_f];\mathbb{R}\right),$ $\alpha(t_0)=x_0$, and $\alpha(t_f)=x_f,$ where $\x_{0,f}=\left(x_{0,f},\mathbf{y}_{0,f}\right)$. Similarly, $\bbeta\in H^1\left([t_0,t_f]; \mathbb{R}^{n-1}\right),$ so that $\bbeta(t_0)=\mathbf{y}_0$ and $\bbeta(t_f)=\mathbf{y}_f$.
	
	\item We denote the Jacobian of a function $\mathbf{F}(x,\mathbf{y})$ as $\nabla \mathbf{F}(x,\mathbf{y})$. For the Jacobian comprised of only the partial derivatives with respect to the components of $\mathbf{y},$ we write $\nabla_{\mathbf{y}}\mathbf{F}(x,\mathbf{y})$.
	
	\item We indicate that a variable is being treated as fixed or pre-determined by separating it with a semicolon. E.g., if the value of $y$ has been set, then $\mathbf{F}(x;\,y)$ is only a function of $x$.
\end{enumerate}

\noindent We also make the following assumptions for each smooth vector field $\mathbf{F}^{\pm}.$ We will use these assumptions to prove two lemmas for the mollified vector field $\fscaled$ in Section \ref{sec:most-probable-paths} which are necessary to establish existence of minimizers to the FW functional. Note, these assumptions also apply in the non-autonomous case but the explicit dependence on time in the argument of $\textbf{F}$ is suppressed. 

\begin{assumptions*}
	Let $\mathbf{F}^{\pm} : \R^n\setminus \{x=0\} \mapsto \R^n$ be smooth vector fields. We assume the following properties about $\mathbf{F}^{\pm}$.
	\begin{enumerate}
		\item \label{ass:growth-conditions}\textit{(Growth conditions)} There exist $R_1,c_1,c_3,c_5>0$ and $c_2,c_4,c_6\in \mathbb{R}$ such that $c_5<c_1$ and for some $1<p<\infty$,  $|\x |>R_1$ implies
		\begin{equation}
		\begin{cases}
		| \mathbf{F}^{-}(\x) | \geq c_1 |\x |^p + c_2, & \text{if } x< 0,\\
		| \mathbf{F}^{+}(\x) | \geq c_1 |\x |^p + c_2, &\text{if } x> 0,
		\end{cases} \label{ass:poly-growth}
		\end{equation}
		and 
		\begin{equation}
		\begin{cases}
		\left| \frac{\partial \mathbf{F}^{-}}{\partial x}(\x) \right| \leq c_3 |\x |^p + c_4, & \text{if } x< 0,\\
		\left| \frac{\partial \mathbf{F}^{+}}{\partial x}(\x) \right| \leq c_3 |\x |^p + c_4, &\text{if } x> 0,
		\end{cases}\label{ass:derivative-decay}
		\end{equation}
		and 
		\begin{equation}
		|F^+(0,\mathbf{y})-F^{-}(0,\mathbf{y})|\leq c_5 |\mathbf{y}|^p+c_6. \label{ass:jump}
		\end{equation}
		
		\item \label{ass:inward-flow} \textit{(Asymptotically inward-flowing)} 
		There exist $R_2, c_7>0$ such that $|\x|>R_2$ implies $\mathbf{F}^{\pm}(\x)\neq 0$ and 
		\begin{equation}
		\begin{cases}
		\langle \mathbf{F}^{-}(\x),\mathbf{r(\x})\rangle <-c_7 |\mathbf{F}^{-}(\x)|, & \text{if } x\leq 0,\\
		\langle \mathbf{F}^{+}(\x),\mathbf{r}(\x)\rangle <-c_7|\mathbf{F}^{+}(\x)|, & \text{if } x\geq 0,
		\end{cases} \label{ass:flow}
		\end{equation}
		where $\mathbf{r}(\x)=\x/|\x|$ is the normalized outward-pointing radial vector at $\x$. 
	\end{enumerate}
\end{assumptions*}

\section{Mollified Freidlin-Wentzell rate functional}
\label{sec:most-probable-paths}
As discussed in Section \ref{sec:intro}, the standard form of the Freidlin-Wentzell rate functional cannot be directly used to calculate most probable paths in piecewise-smooth systems. In this section, we seek to ameliorate this deficiency by studying the convergence of the minimizer of the Freidlin-Wentzell functional, where we approximate Equation \eqref{eq:Intro:PWVF} as a smooth vector field obtained by mollification in the $x-$direction by a smooth, compactly supported, symmetric kernel. Specifically, by smoothing $\mathbf{F}$ via mollification in $x$ over a region of width $2\e$, we obtain a sequence of functionals $\IsigMol$ and a corresponding sequence of minimizers $\balpha_{\e}\in \A$, which we prove converge --- up to a subsequence --- to a minimizer of a limiting functional $\Ilim$ in the limit as  $\e\rightarrow 0$.

\subsection{Mollified functional and existence of a minimum}
\label{sec:min-exists}

We begin this section by recalling the definition of a (Friedrichs) mollifier, following the presentation in \cite{evans2010partial}. Let $\zeta:\R\mapsto\R$ be a smooth even bump function, defined by
\begin{equation}\label{eq:mollifier}
\zeta(x) = \begin{cases}
Ch(x) , & |x| < 1, \\
0, & |x| \geq 1,
\end{cases}    
\end{equation}
where $h:\mathbb{R}\rightarrow \mathbb{R}$ is a smooth function satisfying $h(x)\geq 0$, $h(-1)=h(1)=0$ and all of its derivatives vanish at $x=-1$ and $x=1$, e.g., $h(x)=\exp\left(1/\left( |x|^2 - 1 \right)\right)$. Additionally, $C>0$ is determined so that  $\int_{-\infty}^{\infty}\zeta(x) \, dx = 1$; i.e., $C=1/\int_{-\infty}^{\infty}\,h(x)dx$.  For each $\e>0,$ we define a sequence of functions
$\zeta_{\e}(x)= \zeta\left(x/\e\right)/\e$ which satisfy $\int_{-\e}^{\e} \zeta_{\e}(x) \, dx = 1$ and are compactly supported on $[-\e,\e]$. The \textit{mollification} $\textbf{F}^{\e}$ of $\textbf{F}$ in $x$ is then defined by convolution with $ \zeta_{\e}$ in the $x-$direction: 
\begin{equation}
\begin{aligned}
\textbf{F}^{\e}(x,\mathbf{y}) =  \zeta_{\e}(x)*\textbf{F}(x,\mathbf{y}) 
=& \displaystyle \int_{-\infty}^{0}  \zeta_{\e}(x-s)\mathbf{F}^-(s,\mathbf{y}) \, ds + \int_{0}^{\infty} \zeta_{\e}(x-s)\mathbf{F}^+(s,\mathbf{y}) \, ds, \\
=& \begin{cases}
\displaystyle \int_{-\e}^{\e} \zeta_{\e}(u)\mathbf{F}^-(x-u,\mathbf{y}) \, du , &x\leq-\e, \\
\displaystyle \int_{-\e}^{x} \zeta_{\e}(u)\mathbf{F}^+(x-u,\mathbf{y}) \, du + \int_{x}^{\e}  \zeta_{\e}(u)\mathbf{F}^-(x-u,\mathbf{y}) \, du , &|x|<\e, \\
\displaystyle \int_{-\e}^{\e} \zeta_{\e}(u)\mathbf{F}^+(x-u,\mathbf{y}) \, du , &x\geq\e, \\
\end{cases}\\
=&\int_{-\e}^{\e}\zeta_{\e}(u)\left(\mathbf{F}^{-}(x-u,\mathbf{y})\mathbbm{1}_{\{u\geq x\}}(u)+\mathbf{F}^{+}(x-u,\mathbf{y})\mathbbm{1}_{\{u\leq x\}}(u)\right)du, \label{eq:Feps} 
\end{aligned}
\end{equation}
where $\mathbbm{1}_{A}$ is the indicator function on the set $A$. The \emph{mollified  Freidlin-Wentzell rate functional}  $\IsigMol:\A\mapsto \mathbb{R}$ is then defined by:
\begin{equation}\label{eq:FWMol}
\IsigMol[\balpha]=\int_{t_0}^{t_f}\left| \dot{\balpha}-\fscaled(\balpha)\right|^2\,dt.  
\end{equation}

Note, by construction, $\fscaled(\x)$ is a smooth function that converges pointwise to $\mathbf{F}(\x)$ as $\e\rightarrow 0$ except at $x=0$ \cite{evans2010partial}. However, since $\mathbf{F}(\x)$ is not continuous this convergence is necessarily not uniform and thus it is not \textit{a priori} clear which properties of $\mathbf{F}$ are preserved under mollification. The following lemmas ensure that $\fscaled$ also has $p$-growth and is asymptotically inward-flowing. These properties are crucial to establishing the existence of a minimum for $\IsigMol,$ but 
to streamline our main results, the proofs of these lemmas are in Appendix A.

\begin{lemma}[$p-$growth]\label{lem:p-growth} 
	There exists $\e^*$ such that for all $\e>0$ satisfying $\e<\e^*$, there exist $R^{\e}_1, c_1^{\e}>0$ and $c_2^{\e}\in \mathbb{R}$ such that $|\x|>R_1^{\e}$ implies 
	\[
	|\fscaled(\x)|\geq c_1^{\e}|\x|^p+c_2^{\e}
	\]
	for some $p>1$.
\end{lemma}

\begin{lemma}[Asymptotically inward-flowing]\label{lem:inward-Feps} 
	For all $\e>0,$ there exist $R_2^{\e},c_7^{\e}>0$ such that $|\x|>R_2^{\e}$ implies $\fscaled(\x)\neq 0$ and
	\[
	\langle \fscaled(\x),\mathbf{r}(\x)\rangle <-c_7^{\e}|\fscaled(\x)|,
	\]
	where $\mathbf{r}(\x)=\x/|\x|$ is the normalized outward-pointing radial vector at $\x$.
\end{lemma}

We now prove the existence of a minimizer of $\IsigMol$ using Tonelli's direct method of the calculus of variations; see, e.g.,  \cite{jost1998calculus, braides2002gamma}. This procedure consists of first proving that a minimizing sequence of $\IsigMol$ is bounded with respect to the $H^1$ norm  
and thus has a weakly convergent subsequence in the $H^1$ topology. Note, it is necessary to work in the $H^1$ topology since closed and bounded sets are not necessarily compact in the strong topology. Nevertheless, upon passing to a subsequence, we can establish a limit point for the minimizing sequence. Finally, by proving that  $\IsigMol$ is lower semi-continuous with respect to weak convergence in the $H^1$ topology, it follows that limit points of the minimizing sequence are indeed minimizers.    

To begin the procedure for the direct method of the calculus of variations outlined in the previous paragraph, we first recall the definition of weak convergence in the $H^1$ topology in the context of our problem. 

\begin{definition}[Weak convergence in $H^1$ \cite{evans1990weak}]
	A sequence $\{\balpha_m\}_{m=1}^{\infty} \subset H^1([t_0,t_f];\mathbb{R}^n)$ converges weakly to $\balpha^*\in H^1([t_0,t_f];\mathbb{R}^n)$, written 
	\[
	\balpha_m \rightharpoonup \balpha^* \;\; \mbox{in} \;\; H^1([t_0,t_f];\mathbb{R}^n) \qquad \mbox{or} \qquad \balpha_m \wlim \balpha^*,
	\]
	provided $\balpha_m \rightharpoonup \balpha^*$ in $L^2([t_0,t_f];\mathbb{R}^n)$ and $\dot{\balpha}_m \rightharpoonup \dot{\balpha}^*$ in $L^2([t_0,t_f];\mathbb{R}^n)$. That is, for all $\bm{\omega}\in L^2([t_0,t_f];\mathbb{R}^n)$,
	\[
	\lim_{m\rightarrow  \infty}\int_{t_0}^{t_f} \langle \balpha_m, \bm{\omega} \rangle\,dt = \int_{t_0}^{t_f} \langle \balpha^*, \bm{\omega} \rangle \,dt \qquad \mbox{and} \qquad \lim_{m\rightarrow \infty}\int_{t_0}^{t_f} \langle \dot{\balpha}_m, \bm{\omega} \rangle \,dt = \int_{t_0}^{t_f} \langle \dot{\balpha}^*,\bm{\omega} \rangle\,dt. 
	\]
\end{definition}

The following lemma establishes that boundedness of $\IsigMol[\balpha_{m}]$ implies boundedness of $\balpha_{m}$ with respect to the $H^1$ norm.

\begin{lemma}\label{lem:bdAimpliesbdH1} There exists $\e^*>0$ such that if $\e<\e^*$ and if $\balpha_m\in \mathcal{A}$ satisfies $\IsigMol[\balpha_m]<M_1$ for some $M_1>0$, then there exists $M_2>0$ such that $\|\balpha_m\|_{H^1}<M_2$.
\end{lemma}
\begin{proof}
	For contradiction, suppose $\balpha_m$ is unbounded with respect to the $H^1$ norm. Therefore, there exists a subsequence $\balpha_{m_k}$ such that $\|\balpha_{m_k}\|_{H^1}\rightarrow \infty$. 
	
	1. Suppose there exists $M>0$ such that $\|\balpha_{m_k}\|_{\infty}<M$. Consequently, since $\fscaled$ is smooth there exists $C_1>0$ such that $|\fscaled(\balpha_{m_k})|<C_1$. By the Cauchy-Schwarz inequality,
	\[
	\begin{aligned}
	\IsigMol&=\int_{t_0}^{t_f} \left[ |\dot{\balpha}_{m_k}|^2-2 \langle \balpha_{m_k},\fscaled(\balpha_{m_k})\rangle +|\fscaled(\balpha_{m_k})|^2\right]\,dt \\
	&\geq \int_{t_0}^{t_f} \left[ |\dot{\balpha}_{m_k}|^2-2 \langle \balpha_{m_k},\fscaled(\balpha_{m_k})\rangle \right]\,dt\\
	& \geq \|\dot{\balpha}_{m_k}\|^2-2C_1(t_f-t_0)^{\frac{1}{2}}\|\dot{\balpha}_{m_k}\|.
	\end{aligned}
	\]
	Therefore, since $\|\balpha_{m_k}\|_{H^1}\rightarrow \infty$, it follows from Poincare's inequality that $\|\dot{\balpha}_{m_k}\|\rightarrow \infty$. Thus, $\IsigMol\rightarrow \infty$, which is a contradiction. 
	
	2. Suppose $\|\balpha_{m_k}\|_{\infty}$ is unbounded. Upon passing to another subsequence which we also label $\balpha_{m_k}$, it follows that $\|\balpha_{m_k}\|_{\infty}\rightarrow \infty$. Applying the elementary inequality $a^2+b^2\geq 2ab$, it follows that
	\[
	\begin{aligned}
	\IsigMol[\balpha_{m_k}] &\geq 2\int_{t_0}^{t_f}\left[|\dot{\balpha}_{m_k}| | \fscaled(\balpha_{m_k})|- \langle \balpha_m, \fscaled(\balpha_{m_k})\rangle\right]dt \\
	&= 2\int_{t_0}^{t_f}\left[|\dot{\balpha}_{m_k}| | \fscaled(\balpha_{m_k})|\sin^2\left(\frac{\theta}{2}\right)\right]dt,
	\end{aligned}
	\]
	where $\theta(t)$ is the angle between $\dot{\balpha}_{m_k}(t)$ and $\fscaled(\balpha_{m_k}(t))$. 
	
	We now choose $\e^*,R_1^{\e}$ as in Lemma \ref{lem:p-growth} and $R_2^{\e}$ as in Lemma \ref{lem:inward-Feps}, set $R^{\e}=\max\{R_1^{\e},R_2^{\e}\}$, and let $A_k=\{t: |\balpha_{m_k}(t)|>R^{\e}$ and $\frac{d}{dt}|\balpha_{m_k}(t)|\geq 0\}$. It follows from Lemmas \ref{lem:p-growth} and \ref{lem:inward-Feps} that there exists a constant $C_2>0$ such that
	\[
	\begin{aligned}
	\IsigMol[\balpha_{m_k}] &\geq C_2 \int_{A_k} |\dot{\balpha}_{m_k}| |\fscaled(\balpha_{m_k})|dt \\
	&\geq C_2 \int_{A_k}|\dot{\balpha}_{m_k}| |\balpha_{m_k}|^p\,dt \\
	&= C_2 \int_{A_k}|\dot{\balpha}_{m_k}| |\balpha_{m_k}| |\balpha_{m_k}|^{p-1}\,dt.
	\end{aligned}
	\]
	Applying the Cauchy-Schwarz inequality and the fundamental theorem of calculus, we obtain 
	\[
	\begin{aligned}
	\IsigMol[\balpha_{m_k}] &\geq C_2 \int_{A_k}\langle \dot{\balpha}_{m_k}, \balpha_{m_k}\rangle |\balpha|^{p-1}dt\\
	&= C_2\int_{A_k} \frac{d}{dt} |\balpha|^{p+1}dt\\
	&= C_2\left(\max_{t\in [t_0,t_f]}|\balpha_{m_k}|^{p+1}-(R^{\e})^{p+1}\right).
	\end{aligned}
	\]
	Since $\|\balpha_{m_k}\|_{\infty}\rightarrow \infty$ implies $\max_{t\in [t_0,t_f]}|\balpha_{m_k}(t)|\rightarrow \infty$, it follows that $\IsigMol[\balpha_{m_k}]\rightarrow \infty$, which is a contradiction.
	It follows from items 1 and 2 that $\balpha_m$ is bounded with respect to the $H^1$ norm.
\end{proof}

Finally, we use the direct method of the calculus of variations to prove the existence of a minimum. 
\begin{theorem}[Existence of minimizer of the mollified functional]\label{lem:min-seq-inf}There exists an $\e^*>0$ such that if $\e<\e^*$, then there exists an ${\balpha}^*_{\e}\in \A$ such that for all $\balpha\in \A$,
	\[
	\IsigMol[{\balpha}^*_{\e}]\leq \IsigMol [\balpha].
	\]
\end{theorem}

\begin{proof}
	Let $\e>0,$ and let $\balpha_m$ be a minimizing sequence of $\IsigMol[\balpha]$; i.e.,
	\[
	\lim_{m\rightarrow \infty} \IsigMol[\balpha_m] = \inf_{\balpha\in \mathcal{A}^{(t_0,t_f)}} \IsigMol[\balpha].
	\]
	Consequently, there exists $M_1>0$ such that $\IsigMol[\balpha_m]<M_1$ and thus by Lemma \ref{lem:bdAimpliesbdH1} there exists $M_2>0$ such that $\|\balpha_m\|_{H^1}<M_2$. Thus, by the Banach–Alaoglu theorem and the Rellich-Kondrachov theorem  there exists $\balpha^*_{\e}\in \A$ and a subsequence $\balpha_{m_k}$ such that $\dot{\balpha}_{m_k}\overset{L^2}{\rightharpoonup} \dot{\balpha}^*_{\e}$ and $\balpha_{m_k} \overset{L^\infty}{\rightarrow} {\balpha}^*_{\e}$ as $k\rightarrow\infty$ \cite{adams2003sobolev}. Since the integrand appearing in $\IsigMol$ is convex with respect to $\dot{\balpha},$ it follows that $\IsigMol$ is lower semicontinuous with respect to weak convergence \cite{evans1990weak}, and therefore
	\[
	\IsigMol[{\balpha}^*_{\e}]\leq \liminf_{k\rightarrow \infty}\IsigMol[\balpha_{m_k}]=\lim_{m\rightarrow \infty}\IsigMol[\balpha_m]=\inf_{\balpha\in \A}\IsigMol[\balpha].
	\]
\end{proof}

\subsection{Euler-Lagrange equations}
\label{sec:EL-equations}

Now that the existence of minimizers is established for the Freidlin-Wentzell rate $\IsigMol[\balpha_{\e}]$ corresponding to the mollified system $\dot{\x}=\fscaled$, we can establish a necessary condition for minimizers, which is that $\balpha\in \ASig$ is a $C^2$ integral curve of the flow generated by the following system of Euler-Lagrange (EL) equations:
\begin{equation}\label{eq:EL}
\begin{aligned}
\ddot{\x} = \mathbf{F}^{\e}_t + \left(\nabla\fscaled - \nabla\fscaled^{\intercal}\right)\dot{\x}+\nabla\fscaled^{\intercal}\fscaled,
\end{aligned}
\end{equation}
where the subscript $t$ denotes the partial derivative with respect to time. It is convenient to introduce the scaling $z=x/\e,$ $z\in[-1,1],$ so the EL Equations \eqref{eq:EL} become the fast-slow system
\[
\begin{aligned}
\e\ddot{z} &= F^{\e}_{1,t} + \left(\nabla_{\mathbf{y}}F^{\e}_{1} - \mathbf{G}^{\e}_z\right)\mathbf{\dot{y}}+ \mathbf{G}^{\e\,\intercal}_z \mathbf{G}^{\e}, \\
\ddot{\mathbf{y}} &= \mathbf{G}^{\e}_{t} + \left( \mathbf{G}^{\e}_z - \nabla_{\mathbf{y}}F^{\e}_{1} \right)\e\dot{z} + \left( \nabla_{\mathbf{y}}\mathbf{G}^{\e} - \nabla_{\mathbf{y}}\mathbf{G}^{\e\,\intercal} \right)\dot{\mathbf{y}} + \nabla_{\mathbf{y}}\mathbf{G}^{\e\,\intercal} \mathbf{G}^{\e}.
\end{aligned}
\]
Introducing the conjugate momenta $\varphi=\e\dot{z}-F^{\e}_1$ and $\bpsi = \dot{\mathbf{y}}-\mathbf{G}^{\e}$ leads to the Hamiltonian form,
\begin{equation}\label{eq:EL-rescaled-Hamiltonian}
\begin{aligned}
\e\dot{z} &= F^{\e}_1 + \varphi, \\
\dot{\mathbf{y}} &=\mathbf{G}^{\e} + \mathbf{\bpsi}, \\
\dot{\varphi} &= -\varphi F^{\e}_{1,z} - \langle \mathbf{\bpsi}, \mathbf{G}^{\e}_z \rangle, \\
\mathbf{\dot{\bpsi}} &= -\varphi\nabla_{\mathbf{y}}F^{\e}_1 - \langle \mathbf{\bpsi}, \nabla_{\mathbf{y}}\mathbf{G}^{\e} \rangle.
\end{aligned}
\end{equation}
The most probable path of System \eqref{eq:Intro:SDE},\eqref{eq:Intro:PWVF} corresponds to solutions of the EL Equations \eqref{eq:EL-rescaled-Hamiltonian} subject to boundary conditions for appropriate values of $t$. In the case studies of Sections \ref{sec:case-study-2D} and \ref{sec:case-study-non-autonomous}, we use the gradient flow to numerically solve System \eqref{eq:EL-rescaled-Hamiltonian} when an explicit solution is not available.

\subsection{Limiting functional}
\label{sec:limiting-functional}

In this section, we consider the problem of computing an appropriate limiting functional $\Ilim:\A \mapsto \mathbb{R}$. Theorem  \ref{lem:min-seq-inf} guarantees the existence of a minimizer of $\IsigMol[\balpha]$ for fixed $\e$ and is  given by
$
{\balpha}^*_{\e} = \argmin_{\balpha\in\A} \IsigMol[\balpha(t)].
$
If we now consider the sequence of minimizers ${\balpha}^*_{\e}$ and suppose there exists $\balpha^*\in \A$ such that $\balpha^*_{\e}\wlim \balpha^*$, an appropriate limiting functional $\Ilim$ should have the property that $\balpha^*$ minimizes $\Ilim$. Introduced by de Giorgi, the $\Gamma$-limit is the precise notion of a limiting functional that ensures this property \cite{braides2002gamma, dal2012introduction}.  The following definition of a $\Gamma$-limit is adapted from \cite{braides2002gamma} for our specific problem.
\begin{definition}[$\Gamma-$Convergence \cite{braides2002gamma}]
	We say that $\IsigMol$ $\Gamma-${\upshape converges } to $\Ilim:\A\mapsto \overline{\R}$ with respect to $H^1$ weak convergence if for all $\balpha\in \A$, we have
	\begin{enumerate}
		\item (liminf inequality) for every sequence $\balpha_{\e}\in \A$ satisfying $\balpha_{\e}\overset{H^1}{\rightharpoonup} \balpha$, 
		\[
		\Ilim[\balpha]\leq \liminf_{\e \rightarrow 0}\IsigMol[\balpha_{\e}];
		\]
		\item (recovery sequence) there exists a sequence $\balpha_{\e}\overset{H^1}{\rightharpoonup}\balpha$ such that
		\[
		\lim_{\e \rightarrow 0}\IsigMol[\balpha_{\e}]=\Ilim[\balpha].
		\]
	\end{enumerate}
	If $\IsigMol$ $\Gamma$-converges to $\Ilim$ we write $\displaystyle \Gamma-\lim_{\e \rightarrow 0}\IsigMol=\Ilim$.
\end{definition}

Showing both the liminf inequality and the existence of a recovery sequence establish that both the limsup and the liminf of $\IsigMol$ over all weakly converergent sequences exists and is equal to $\Ilim$. Indeed, these two conditions ensure that $\IsigMol$ has a common lower bound, the lower bound is optimal, and the limiting functional is precisely this lower bound.   Moreover, the definition of $\Gamma$-convergence is constructed in such a manner that ensures the convergence of minimizers of $\IsigMol$ to the minimizer of $\Ilim$ if an additional equicoecervity condition is satisfied. The following theorem adapted from \cite{braides2002gamma} ensures this convergence for our specfic problem.
\begin{theorem}\label{thm:min-exists} 
	Suppose $\IsigMol$ $\Gamma$-converges to $\Ilim:\A \mapsto \overline{\R}$ as $\e\rightarrow 0$ with respect to $H^1$ weak convergence and suppose for all $\balpha^*_{\e}$ that minimize $\IsigMol$ there exists a constant $M>0$ such that $\|\balpha_{\e}^*\|_{H^1}<M$. Then there exists an $\balpha^*\in \A$ such that
	\[
	\inf_{\balpha \in \A} \Ilim[\balpha]=\Ilim[\balpha^*].
	\]
	Moreover, every limit of a subsequence of $\balpha^*_{\e}$ as $\e\rightarrow 0$ is a minimizer of $\Ilim$.
\end{theorem}

One of the challenges with proving a $\Gamma$-limit is that the target limiting functional must be formed from an educated guess based on analysis of minimizers. In the system of interest, it is natural to expect that the appropriate $\Gamma$-limit will correspond to the standard FW functional except for intervals of time on which the curve tracks $\Sigma$. However, for $\balpha\in \A$ satisfying $\{t\in \mathbb{R}: \alpha(t)=0\}$ has nonzero measure, the recovery sequence must be constructed to minimize the functional for times in which $\balpha(t) \in \Sigma^{\e}=\{\mathbf{x}\in \mathbb{R}^n: |x|< \e\}$. To compute the appropriate $\Gamma$-limit for $\IsigMol,$ we introduce the following intervals of time:
\[
\begin{split}
\I[\balpha]&=\{t\in (t_0,t_f): \balpha(t)\not\in\Sigma\}=\{t\in (t_0,t_f): \alpha_1(t)\neq 0\},\\
\I^{\e}[\balpha]&=\{t\in (t_0,t_f): \balpha(t)\not\in\Sigma^{\e}\}=\{t\in (t_0,t_f):|\alpha_1(t)|\geq \e \},\\
\I_{\Sigma}[\balpha]&=\{t\in (t_0,t_f): \balpha(t)\in\Sigma\} = \{t\in (t_0,t_f): \alpha_1(t) =0\},\\
\I_{\Sigma_{+}}[\balpha]&=\{t\in (t_0,t_f): \balpha(t)\in\Sigma_{+}\} = \{t\in (t_0,t_f): \alpha_1(t) =0 \mbox{ and } F^{\pm}_1(0,\bbeta)>0\},\\
\I_{\Sigma_{-}}[\balpha]&=\{t\in (t_0,t_f): \balpha(t)\in\Sigma_{-}\} = \{t\in (t_0,t_f): \alpha_1(t) =0 \mbox{ and } F^{\pm}_1(0,\bbeta)<0\},\\
\I^{\e}_{\Sigma}[\balpha]&=\{t\in (t_0,t_f): \balpha(t)\in\Sigma^{\e}\} = \{t\in (t_0,t_f): |\alpha_1(t)|<\e\},\\
\I_{\pm}[\balpha]&=\{t\in (t_0,t_f): \balpha(t)\in S_{\pm}\},\\
\I_{\pm}^{\e}[\balpha]&=\{t\in (t_0,t_f): \balpha(t)\in S_{\pm}^{\e}\},\\
\end{split}
\]
as well as introduce the function $\Lambda:\R\mapsto\R$, defined by
\begin{equation}\label{eq:Lambda}
\Lambda(z)=\frac{1}{2}+\int_0^z \zeta(u)\,du.
\end{equation}
Note that $\Lambda(z)=1$ if $z>1,$ $\Lambda(z)=0$ if $z<-1,$ and $\Lambda(z)$ smoothly transitions from $0$ to $1$ on the interval $[-1,1]$. We now state the following $\Gamma-$limit result.

\begin{theorem}\label{thm:gamma-convergence} 
	There exists $\Ilim :\A\mapsto \overline{\R}$ such that
	\[
	\Gamma-\lim_{\e \rightarrow 0}\IsigMol=\Ilim,
	\]
	where for $\balpha=(\alpha,\bbeta)\in \A$,
	\begin{equation}\label{eq:gamma-limit}
	\begin{aligned}
	\Ilim[\balpha]=&\int_{\I[\balpha]}\|\dot{\balpha}(t)-\mathbf{F}(\balpha(t))\|^2\,dt \\
	&+\int_{\I_{\Sigma}[\balpha]}\min_{\lambda\in [0,1]}\left\{
	\left[ \lambda F_1^+(0,\bbeta) + (1-\lambda)F_1^-(0,\bbeta) \right]^2 \right. \\
	&\left. \mathmakebox[\widthof{$ +\int_{\I_{\Sigma}[\balpha]}\min\qquad$ }][l]{}  + \left| \dot{\bbeta} - \lambda\mathbf{G}^+(0,\bbeta) - (1-\lambda)\mathbf{G}^-(0,\bbeta) \right|^2
	\right\}\,dt.
	\end{aligned}
	\end{equation}
\end{theorem}

\begin{proof} 
	For fixed $\balpha =(\alpha,\bbeta)\in \A$, we first show the existence of a recovery sequence, $\balpha_{\e}$. We will construct $\balpha_{\e}$ over three intervals, $\I^{\e}_+[\balpha],$ $\I^{\e}_-[\balpha],$ and $\I^{\e}_{\Sigma}[\balpha],$ by
	\[
	\begin{aligned}
	\balpha_{\e}(t) &= \balpha(t)\mathbbm{1}_{\{t\in \I^{\e}_+[\balpha]\}} + \balpha(t)\mathbbm{1}_{\{t\in \I^{\e}_-[\balpha]\}} + (\e z(t),\bbeta(t))\mathbbm{1}_{\{t\in \I_{\Sigma}^{\e}[\balpha]\}} \\
	&= \balpha(t)\mathbbm{1}_{\{\balpha\in S^{\e}_+\}} + \balpha(t)\mathbbm{1}_{\{\balpha\in S^{\e}_-\}} + (\e z(t),\bbeta(t))\mathbbm{1}_{\{(\e z,\bbeta)\in \Sigma^{\e}\}}, 
	\end{aligned}
	\]
	where $z(t)$ is given by 
	\begin{equation}\label{eq:z-def}
	\begin{aligned}
	z = \argmin_{z'\in [-1,1]} &\left\{ 
	\left[ \Lambda(z')F_1^+(0,\bbeta) + \left(1-\Lambda(z')\right) F_1^-(0,\bbeta) \right]^2 \right.\\
	&\mathmakebox[\widthof{$ \left\{ 
		\left[\right.\right.d$ }][l]{} + \left. \left| \dot{\bbeta} - \Lambda(z')\mathbf{G}^+(0,\bbeta) - \left(1-\Lambda(z')\right)\mathbf{G}^-(0,\bbeta) \right|^2    
	\right\}.
	\end{aligned}
	\end{equation}
	
	To prove convergence of $\IsigMol[\balpha_{\e}]$ to $\Ilim[\balpha]$ and $\balpha_{\e}\wlim \balpha,$ we first consider the interval of time $\mathcal{I}_+^{\varepsilon}[\balpha]$. Since $\balpha_{\e}(t)\mathbbm{1}_{\{\balpha_{\e}\in S^{\e}_+\}}=\balpha(t)\mathbbm{1}_{\{\balpha\in S^{\e}_+\}}$ it follows that $\balpha_{\e}(t)\mathbbm{1}_{\{\balpha_{\e}\in S^{\e}_+\}}\wlim\balpha(t)\mathbbm{1}_{\{\balpha\in S^{\e}_+\}}$. Moreover, we have
	\[
	\begin{aligned}
	\int_{\I^{\e}_+[\balpha_{\e}]} & \|\dot{\balpha}_{\e}-\fscaled(\balpha_{\e})\|^2\,dt - \int_{\I_+[\balpha]}\|\dot{\balpha}-\mathbf{F}^+(\balpha)\|^2\,dt \\
	=& \int_{\I^{\e}_+[\balpha] \cap \I_+[\balpha]} \left( 2\langle \dot{\balpha},\mathbf{F}^+(\balpha)-\fscaled (\balpha)\rangle + |\fscaled(\balpha)|^2 - |\mathbf{F}^+(\balpha)|^2 \right)\,dt \\
	&+ \int_{\I^{\e}_+[\balpha] \setminus \I_+[\balpha]} \|\dot{\balpha}-\fscaled(\balpha)\|^2\,dt
	- \int_{\I_+[\balpha] \setminus \I^{\e}_+[\balpha]} \|\dot{\balpha}-\mathbf{F}^+(\balpha)\|^2\,dt.
	\end{aligned}
	\]
	Since $\fscaled(x,\mathbf{y})\mathbbm{1}_{\{x\geq\e\}} = \int_{-1}^1 \zeta(v)\mathbf{F}^+(x-\e v,\mathbf{y})\,dv\mathbbm{1}_{\{x\geq\e\}}$ and $\int_{-1}^1 \zeta(v)\mathbf{F}^+(x-\e v,\mathbf{y})\,dv$ converges uniformly to $  \mathbf{F}^+(x,\mathbf{y})$ as $\e\rightarrow 0,$ the integral $\int_{\I^{\e}_+[\balpha] \cap \I_+[\balpha]} \langle \dot{\balpha},\mathbf{F}^+(\balpha)-\fscaled(\balpha) \rangle \,dt$ vanishes as $\e\rightarrow 0$. By Minkowski's inequality and the triangle inequality, 
	\[
	\begin{aligned}
	\left| \int_{\I^{\e}_+[\balpha] \cap \I_+[\balpha]} |\fscaled(\balpha)|^2 - |\mathbf{F}^+(\balpha)|^2 \,dt \right| 
	&\leq \int_{\I^{\e}_+[\balpha] \cap \I_+[\balpha]} \left| \|\fscaled(\balpha)\|^2 - \|\mathbf{F}^+(\balpha)\|^2 \right| \,dt\\
	&\leq \int_{\I^{\e}_+[\balpha] \cap \I_+[\balpha]}  \|\fscaled(\balpha) - \mathbf{F}^+(\balpha)\|^2 \,dt,
	\end{aligned}
	\]
	and thus by the same uniform convergence argument as above vanishes as $\e\rightarrow 0$.  Since the measures of $\I^{\e}_+[\balpha] \setminus \I_+[\balpha]$ and $\I_+[\balpha] \setminus \I^{\e}_+[\balpha]$ vanish as $\e \rightarrow 0$ it follows that the remaining integrals also vanish. The analogous case for $\I_-[\balpha]$ is identical. Hence,
	\[
	\begin{aligned}
	\lim_{\e\rightarrow 0}\IsigMol[\balpha_{\e}\mathbbm{1}_{\{\balpha_{\e}\in S^{\e}_{\pm}\}}] &= \lim_{\e\rightarrow 0}\int_{\I^{\e}[\balpha]} \|\dot{\balpha}-\fscaled(\balpha)\|^2\,dt \\
	&= \int_{\I[\balpha]}\|\dot{\balpha}-\mathbf{F}(\balpha)\|^2\,dt \\
	&= \Ilim[\balpha\mathbbm{1}_{\{\balpha\not\in \Sigma\}}].
	\end{aligned}
	\]
	
	We now consider the interval of time $\I^{\e}_{\Sigma}[\balpha]$. Note that using the substitution $v=u/\e,$ we have
	\[
	\begin{aligned}
	\lim_{\e\rightarrow 0} \fscaled(\e z,\mathbf{y})\mathbbm{1}_{\{ |\e z|<\e \}} &= \lim_{\e\rightarrow 0}  \left[ \int_{-1}^{ z}\zeta(v)\mathbf{F}^+(\e z - \e v,\mathbf{y})\,dv + \int_z^1 \zeta(v)\mathbf{F}^-(\e z-\e v,\mathbf{y})\,dv\right]  \\
	&= \Lambda(z)\mathbf{F}^+(0,\mathbf{y}) + \left(1-\Lambda(z)\right)\mathbf{F}^-(0,\mathbf{y}).
	\end{aligned}
	\]
	Let 
	\[
	\mathbf{H}(\mathbf{y};\,z) = \Lambda(z)\mathbf{F}^+(0,\mathbf{y}) + \left(1-\Lambda(z)\right)\mathbf{F}^-(0,\mathbf{y}) = \lambda\mathbf{F}^+(0,\mathbf{y}) + \left(1-\lambda\right)\mathbf{F}^-(0,\mathbf{y}),
	\]
	where we set $z(t)$ as in Equation \eqref{eq:z-def}, so that $\lambda = \Lambda(z)\in[0,1]$ is a constant.
	Thus, in general $\lim_{\e\rightarrow 0} \fscaled(\e z,\mathbf{y}) = \mathbf{H}(z,\mathbf{y})$, where $z$ varies in $[-1,1],$ but we evaluate $\mathbf{H}$ at a particular $z-$value to correspond with the $\Gamma-$limit. Since $\balpha_{\e}(t)\mathbbm{1}_{\{\balpha_{\e}\in \Sigma^{\e}\}}=(\e z,\bbeta)\in \A$ it follows that
	\[
	\begin{aligned}
	\int_{\I^{\e}_{\Sigma}[\balpha]} &\|\dot{\balpha}_{\e}-\fscaled(\balpha_{\e})\|^2\,dt - \int_{\I_{\Sigma}[\balpha]}\left\|\dot{\balpha}-\mathbf{H}(\bbeta;\, z)\right\|^2\,dt \\
	=& \int_{\I^{\e}_{\Sigma}\cap \I_{\Sigma}} \|\dot{\balpha}_{\e}-\fscaled(\balpha_{\e})\|^2 - \left\|\dot{\balpha}-\mathbf{H}(\bbeta;\, z)\right\|^2\,dt \\
	& + \int_{\I^{\e}_{\Sigma} \setminus \I_{\Sigma}} \|\dot{\balpha}_{\e}-\fscaled\|^2\,dt
	- \int_{\I_{\Sigma} \setminus \I^{\e}_{\Sigma}} \|\dot{\balpha}-\mathbf{H}\|^2\,dt.
	\end{aligned}
	\]
	As above, the last two integrals vanish as $\e\rightarrow 0,$ since the measure of $\I^{\e}_{\Sigma} \setminus \I_{\Sigma}$ converges to zero and $\I_{\Sigma} \setminus \I^{\e}_{\Sigma}=\varnothing,$ and the integrands are bounded. Since $\I^{\e}_{\Sigma}\cap \I_{\Sigma} = \I_{\Sigma},$ 
	\[
	\begin{aligned}
	\int_{\I^{\e}_{\Sigma}\cap \I_{\Sigma}} &\|\dot{\balpha}_{\e}-\fscaled(\balpha_{\e})\|^2 - \left\|\dot{\balpha}-\mathbf{H}(\bbeta;\, z)\right\|^2\,dt \\ 
	=& \int_{\I_{\Sigma}} |\dot{\balpha}_{\e}|^2-2\langle \dot{\balpha}_{\e}, \fscaled(\balpha_{\e}) \rangle + |\fscaled(\balpha_{\e})|^2 - \left|\dot{\balpha}\right|^2+2\langle \dot{\balpha}, \mathbf{H}(\bbeta;\, z) \rangle -\left|\mathbf{H}(\bbeta;\, z)\right|^2\,dt.
	\end{aligned}
	\]
	Since $\dot{\balpha}_{\e}\mathbbm{1}_{\{\balpha_{\e}\in\Sigma^{\e}\}}\wlim \dot{\balpha}\mathbbm{1}_{\{\balpha\in\Sigma\}}$ and $\fscaled(\balpha_{\e})\mathbbm{1}_{\{\balpha_{\e}\in\Sigma^{\e}\}}\wlim \mathbf{H}(\bbeta;\, z)\mathbbm{1}_{\{\balpha\in\Sigma\}}, $
	\[
	\begin{aligned}
	\int_{\I_{\Sigma}} |\dot{\balpha}_{\e}|^2- |\dot{\balpha}|^2\,dt = \left| \int_{t_0}^{t_f} |\dot{\balpha}_{\e}\mathbbm{1}_{\{\balpha_{\e}\in\Sigma^{\e}\}} |^2 - |\dot{\balpha}\mathbbm{1}_{\{\balpha\in\Sigma\}}|^2 \,dt \right| 
	\leq \int_{t_0}^{t_f}  \|\left(\dot{\balpha}_{\e} - \dot{\balpha}\right)\mathbbm{1}_{\{\balpha_{\e}\in\Sigma^{\e}\}}\|^2 \,dt \rightarrow 0
	\end{aligned}
	\]
	as $\e\rightarrow 0$, and
	\[
	\begin{aligned}
	\lim_{\e\rightarrow 0}\int_{\I_{\Sigma}} |\fscaled(\balpha_{\e})|^2 -\left|\mathbf{H}(\bbeta;\, z)\right|^2\,dt &= \lim_{\e\rightarrow 0}\left| \int_{t_0}^{t_f} |\fscaled\mathbbm{1}_{\{\balpha_{\e}\in\Sigma^{\e}\}} |^2 - |\mathbf{H}\mathbbm{1}_{\{\balpha\in\Sigma\}}|^2 \,dt \right| \\ 
	&\leq \lim_{\e\rightarrow 0}\int_{t_0}^{t_f} \|\left(\fscaled - \mathbf{H}\right)\mathbbm{1}_{\{\balpha_{\e}\in\Sigma^{\e}\}}\|^2 \,dt \\ &= 0,
	\end{aligned}
	\]
	where we have used Minkowski's inequality and the triangle inequality, as above. Also, as a consequence of the uniform boundedness of sequences in $H^1,$ the product of a weakly convergent sequence and a strongly convergent sequence is weakly convergent \cite{evans1990weak}, and thus
	\[
	\begin{aligned}
	\int_{\I_{\Sigma}} 
	\langle \dot{\balpha}, \mathbf{H}(\bbeta;\, z) \rangle
	- \langle \dot{\balpha}_{\e}, \fscaled(\balpha_{\e}) \rangle \,dt = \int_{t_0}^{t_f} 
	\langle \dot{\balpha}, \mathbf{H}(\bbeta;\, z) \rangle \mathbbm{1}_{\{\balpha\in\Sigma\}}
	- \langle \dot{\balpha}_{\e}, \fscaled(\balpha_{\e})\rangle\mathbbm{1}_{\{\balpha_{\e}\in\Sigma^{\e}\}} \,dt \rightarrow 0
	\end{aligned}
	\]
	as $\e\rightarrow 0$. Hence,
	\[
	\begin{aligned}
	\lim_{\e\rightarrow 0}\IsigMol[\balpha_{\e}\mathbbm{1}_{\{\balpha_{\e}\in \Sigma^{\e}\}}] &= \lim_{\e\rightarrow 0} \int_{\I^{\e}_{\Sigma}} \|\dot{\balpha}_{\e}-\fscaled(\balpha_{\e})\|^2\,dt \\ 
	&= \int_{\I_{\Sigma}}\|\dot{\balpha}-\mathbf{H}(\bbeta;z)\|^2\,dt \\
	&= \int_{\I_{\Sigma}}\min_{\lambda\in[0,1]}\left\{ \left|\lambda F^+_1(0,\bbeta) + \left(1-\lambda\right)F^-_1(0,\bbeta)\right|^2 \right. \\ 
	& \left.\mathmakebox[\widthof{$ = \int_{\I_{\Sigma}}\min_{\lambda\in[0,1]}$ }][l]{} +  \left|\dot{\bbeta}-\lambda\mathbf{G}^+(0,\bbeta) - \left(1-\lambda\right)\mathbf{G}^-(0,\bbeta)\right|^2 \right\}\,dt \\
	&= \Ilim[\balpha\mathbbm{1}_{\{\balpha\in {\Sigma}\}}].
	\end{aligned}
	\]
	
	We now prove the liminf inequality. 
	For the portion of $\balpha$ in $S_{\pm},$ the liminf inequality follows directly from the recovery sequence calculation, since we chose $\balpha_{\e}=\balpha$ for any $\balpha\in\A$. For the portion of $\balpha$ in $\Sigma,$ we have that for every $\balpha_{\e},\balpha\in\A$ satisfying $\balpha_{\e}\mathbbm{1}_{\{\balpha_{\e}\in \Sigma^{\e}\}}\wlim\balpha\mathbbm{1}_{\{\balpha\in \Sigma\}},$
	\[
	\begin{aligned}
	\lim_{\e\rightarrow 0} \IsigMol[\balpha_{\e}\mathbbm{1}_{\{\balpha_{\e}\in \Sigma^{\e}\}}] &= \lim_{\e\rightarrow 0} \int_{t_0}^{t_f} |\dot{\balpha}_{\e}\mathbbm{1}_{\{\balpha_{\e}\in \Sigma^{\e}\}} - \fscaled\left(\balpha_{\e}\mathbbm{1}_{\{\balpha_{\e}\in \Sigma^{\e}\}}\right)|^2\,dt \\
	&= \int_{\I_{\Sigma}[\balpha]}|\dot{\balpha}\mathbbm{1}_{\{\balpha\in \Sigma\}}-\mathbf{H}(\balpha\mathbbm{1}_{\{\balpha\in \Sigma\}})|^2\,dt \\
	& \geq \int_{\I_{\Sigma}[\balpha]} \min_{\lambda\in[0,1]}\left\{\left|\lambda F^+_1(0,\bbeta) + \left(1-\lambda\right)F^-_1(0,\bbeta)\right|^2 \right.\\ 
	& \mathmakebox[\widthof{$ = \int_{\I_{\Sigma}[\balpha]} \min\left\{\right.$ }][l]{} +  \left.\left|\dot{\bbeta}-\lambda\mathbf{G}^+(0,\bbeta) - \left(1-\lambda\right)\mathbf{G}^-(0,\bbeta)\right|^2\right\}\,dt.
	\end{aligned}
	\]
	Then by the lower semicontinuity of the rate functional with respect to weak convergence,
	\[
	\begin{aligned}
	\Ilim[\balpha\mathbbm{1}_{\{\balpha\in \Sigma\}}] &= \int_{\I_{\Sigma}} \min_{\lambda\in[0,1]}\left\{\left|\lambda F^+_1(0,\bbeta) + \left(1-\lambda\right)F^-_1(0,\bbeta)\right|^2 \right.\\ 
	& \mathmakebox[\widthof{$ = \int_{\I_{\Sigma}[\balpha]} \min\left\{\right.$ }][l]{} +  \left.\left|\dot{\bbeta}-\lambda\mathbf{G}^+(0,\bbeta) - \left(1-\lambda\right)\mathbf{G}^-(0,\bbeta)\right|^2\right\}\,dt \\ 
	&\leq \liminf_{\e\rightarrow 0} \IsigMol[\balpha_{\e}\mathbbm{1}_{\{\balpha_{\e}\in \Sigma^{\e}\}}].
	\end{aligned}
	\]
	
\end{proof}

In other words, this result shows that in the smooth regions $S_{\pm},$ the minimum of the rate functional $\Isig$ by $\balpha(t),$ where $\balpha$ is the solution to the EL Equations \eqref{eq:EL-rescaled-Hamiltonian}. Then in the switching manifold $\Sigma$, the rate functional is minimized by some $\balpha(t;z)=(0,\bbeta(t;z))\in\A,$ where $z$ is defined using Equation \eqref{eq:z-def}. Note that the critical manifold of the fast-slow System \eqref{eq:EL-rescaled-Hamiltonian} is defined implicitly as solutions in $\mathbb{R}^n$ to following system of equations
\[
M_0 \subseteq \left\{\begin{aligned}
0 &= \Lambda(z)F_1^+(0,\mathbf{y}) + \left(1-\Lambda(z)\right)F_1^-(0,\mathbf{y}) + \varphi \\
0 &= \varphi\left( F_1^+(0,\mathbf{y}) + F_1^-(0,\mathbf{y})\right) + \langle \bpsi, \mathbf{G}^+(0,\mathbf{y}) + \mathbf{G}^-(0,\mathbf{y}) \rangle
\end{aligned}\right\},
\]
so when the conjugate momentum has $\varphi=0,$ $\Lambda(z)F_1^+(0,\mathbf{y}) + \left(1-\Lambda(z)\right)F_1^-(0,\mathbf{y})=0$. Thus, the path that minimizes the rate functional in $\Sigma$ is given by $\balpha=(0,\bbeta),$ where 
\begin{equation}\label{eq:beta-dot}
\dot{\bbeta} = \Lambda(z)\mathbf{G}^+(0,\bbeta) + \left(1-\Lambda(z)\right)\mathbf{G}^-(0,\bbeta).
\end{equation}
Practically, to obtain the most probable path in the switching manifold, we first determine the value of $z$ so that $\Lambda(z)F_1^+(0,\mathbf{y}) + \left(1-\Lambda(z)\right)F_1^-(0,\mathbf{y})=0$ for each $\mathbf{y}\in\R^{n-1}$. This fixes $z$, so we may calculate the path using Equation \eqref{eq:beta-dot}.

To make explicit the connection of this result to the dynamics imposed  on the switching manifold by the Filippov convex combination, recall the definition of sliding flow given by Equations \eqref{eq:sliding-flow},\eqref{eq:gamma}. Here, we may fix $\lambda=\Lambda(z)$ so that the first component of the sliding flow in Equation \eqref{eq:sliding-flow} becomes zero. The remaining components describing the sliding flow in Equation \eqref{eq:sliding-flow} then either coincide with Equation \eqref{eq:beta-dot} or its time-reversed flow, with $t\rightarrow -t.$ That is, Theorem \ref{thm:gamma-convergence} indicates that the intersection of the most probable path with the switching manifold, $\balpha(t)\mathbbm{1}_{\{\balpha\in \Sigma\}},$ may correspond to the solution of the flow given by $\dot{\mathbf{y}}=\mathbf{G}^s(0,\mathbf{y};z),$ where $\mathbf{G}^s$ is the convex combination
\[
\mathbf{G}^s(0,\mathbf{y};z) = \lambda \mathbf{G}^+(0,\mathbf{y};z) + (1-\lambda)\mathbf{G}^-(0,\mathbf{y};z),
\]    
and $\lambda\in [0,1]$. Thus, the most probable path in the switching manifold precisely recovers the Filippov dynamics with this particular form of $\lambda$ (or its time-reversed dynamics). This result is independent of the mollifier imposed, as long as it fits the form defined in Equation \eqref{eq:mollifier}, though the path does depend on the mollifier via $\Lambda(z)$.

Furthermore, this indicates that the functional given by Equation \eqref{eq:gamma-limit} has no contribution during times when the path slides in a sliding region, but it does when the path slides in a crossing region.

\section{Case Study 1: Linear system in 2D}
\label{sec:case-study-2D}

In this section, we explore the application of the most probable path derived in Section \ref{sec:most-probable-paths} to a simple piecewise-linear system in $\R^2$. Consider System \eqref{eq:Intro:SDE},\eqref{eq:Intro:PWVF}, where $\mathbf{F}^{\pm}:\R^2\mapsto \R^2$ are defined as 
\begin{equation}\label{eq:linear-system}
\begin{aligned}
\textbf{F}^+(x,y) &= \begin{pmatrix}
f^+ \\
g^+
\end{pmatrix} = \begin{pmatrix}
a(x - 1) + b(y - \eta) \\
c(x - 1) + y - \eta
\end{pmatrix}, \\
\textbf{F}^-(x,y) &= \begin{pmatrix}
f^- \\
g^-
\end{pmatrix} = \begin{pmatrix}
p(x + 1) + q(y + \eta) \\
r(x + 1) + y + \eta
\end{pmatrix}.
\end{aligned}
\end{equation}
Here, $a,b,c,p,q,r,\eta\in\R$ are parameters. This splits up the plane into two smooth regions separated by the switching manifold $\Sigma=\{(x,y)~|~x=0\}$. This system has been nondimensionalized to scale out a parameter.

This case study is intended to represent a prototypical example of a planar model with a one-dimensional switching manifold and linear or linearizable dynamics near the switch; see, e.g., \cite{berglund2002metastability,arnscheidt2020routes,chen2014large,baule2010stick,piltz2014prey,serdukova2017metastability,monahan2002stabilization}. 
For example, the stochastic Stommel model for thermohaline convection switches between temperature- and salinity-driven ocean circulation \cite{berglund2002metastability,monahan2002stabilization}. 
Here the vector field is nonlinear, but it is reasonable to consider linearizing the vector field near the switch.

We begin by investigating the dynamics of the deterministic skeleton, System \eqref{eq:Intro:DS},\eqref{eq:linear-system}, as the parameter $\eta$ varies. For all values of $\eta,$ there are two stable fixed points separated by the switching manifold $\Sigma$. For most values of $\eta,$ dynamics fall under two broad categories: \textit{(1)} the two fixed points are the only asymptotic states possible for $\x_0\in\R^2\setminus\Sigma$, and \textit{(2)} there is a sliding-cycle surrounding each fixed point, with a crossing-cycle surrounding all sliding-cycles (see Figure \ref{fig:linear-stability}(f)). For intermediate values of $\eta$ we observe \textit{(3)} both standard and nonstandard discontinuity-induced bifurcations. 

Taking the two most prevalent cases \textit{(1)} and \textit{(2)}, we perform Monte Carlo simulations of System \eqref{eq:Intro:SDE},\eqref{eq:linear-system}. In the case of attracting sliding, Monte Carlo simulations match the most probable path derived by minimizing the functional Equation \eqref{eq:gamma-limit}, which slides along $\Sigma$. However, in the case of repelling sliding, Monte Carlo simulations match a local minimizer that does not slide, instead crossing $\Sigma$ at the onset of sliding.

\subsection{Deterministic dynamics}

\begin{figure}[t!]
	\centering
	\begin{subfigure}{.3\textwidth}
		\centering
		\includegraphics[scale=1]{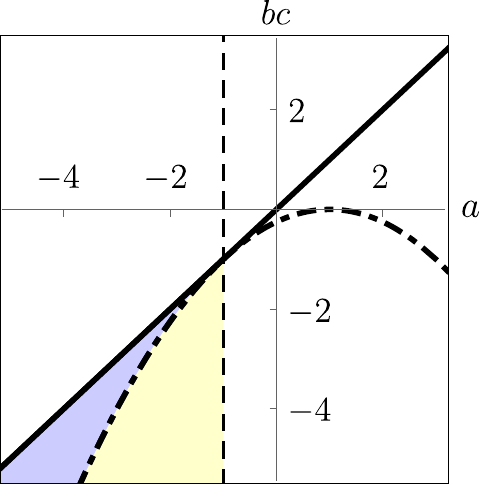}
		\caption{}
	\end{subfigure}%
	\begin{subfigure}{.7\textwidth}
		\centering
		\includegraphics[scale=1]{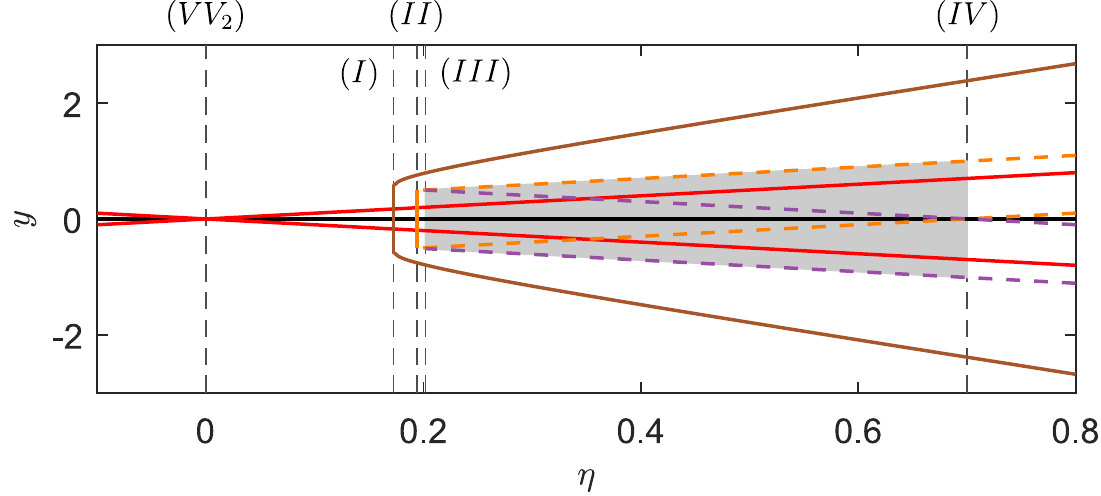}
		\caption{}
	\end{subfigure}\\
	\begin{subfigure}{.25\textwidth}
		\centering
		\includegraphics[scale=1]{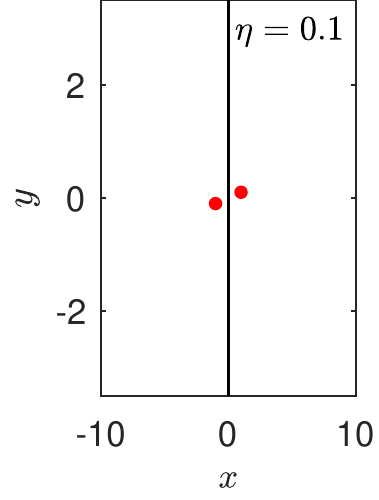}
		\caption{}
	\end{subfigure}%
	\begin{subfigure}{.25\textwidth}
		\centering
		\includegraphics[scale=1]{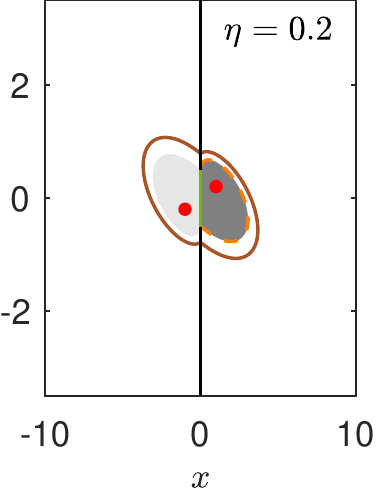}
		\caption{}
	\end{subfigure}%
	\begin{subfigure}{.25\textwidth}
		\centering
		\includegraphics[scale=1]{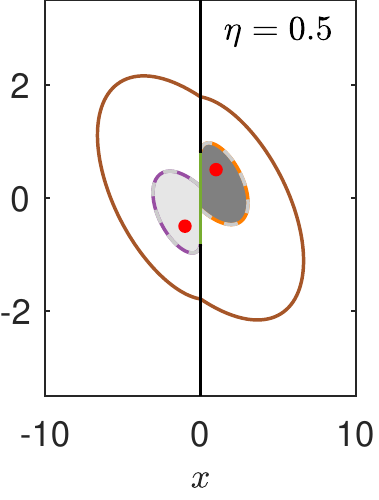}
		\caption{}
	\end{subfigure}%
	\begin{subfigure}{.25\textwidth}
		\centering
		\includegraphics[scale=1]{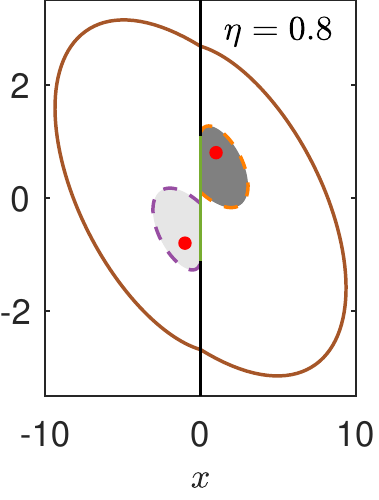}
		\caption{}
	\end{subfigure}
	\caption{(a) Stability plane for $\x_{+}$. Shaded regions have stable nodes (blue) or stable spirals (yellow). (b) A bifurcation diagram for System \eqref{eq:Intro:DS},\eqref{eq:linear-system}. Red lines are stable foci and the black line is a pseudosaddle. Brown curves are the crossing points of the large crossing limit cycle. The orange and purple dashed curves are the upper and lower limits of the sliding portion of the two (unstable) sliding limit cycles. The grey shaded region indicates $\eta-$values where the two sliding limit cycles may have both period-1 and period-2. Vertical dashed lines indicate bifurcations identified: a visible tangency bifurcation $(VV_2)$, a crossing cycle-tangency bifurcation $(I),$ two sliding cycle-tangency bifurcations $(II)$ and $(III),$ and a sliding cycle-nonunique period bifurcation $(IV)$. (c)-(f) Phase portraits for illustrative values of $\eta$. Colors correspond to (b), with the period-2 unstable sliding limit cycle indicated as a dashed gray curve. Parameters used in all plots are $a=p=-2,\, b=q=-7,$ and $c=r=1$.
	}
	\label{fig:linear-stability}
\end{figure}

In System \eqref{eq:Intro:DS},\eqref{eq:linear-system}, there is one equilibrium in each smooth region, at $\x_{\pm}=(\pm 1,\pm\eta)$. Figure \ref{fig:linear-stability}(a) summarizes the linear stability analysis of $\x_{+}$; due to the symmetry of the system, the linear stability of $\x_{-}$ is identical if we consider the vertical axis as $qr$ and the horizontal axis as $p$. The solid line separates saddle ($a<bc$) from non-saddle equilibria. For non-saddle equilibria, the dashed line further separates stable equilibria ($a<-1$) from unstable equilibria. Stable spirals occur in the yellow region and stable nodes occur in the blue region, which are separated by the curve $bc = -(a-1)^2/4$.

Figure \ref{fig:linear-stability}(b) shows the bifurcations in the system as $\eta$ varies, in the case where $\x_{\pm}$ are stable spirals ($\eta$ does not affect the linear stability  of $\x_{\pm}$). The remainder of this section discusses dynamics shown in this bifurcation diagram. We identify bifurcations using existing terminology where possible, though some bifurcations do not appear to have existing classification. Note, although we discuss the time-reversed dynamics in this section, in general most probable paths are not time-reversible unless there is a gradient structure \cite{grafke2019string}.

For large $\eta,$ in addition to the two stable foci $\textbf{x}_{\pm},$ the system has three limit cycles of the piecewise-smooth flow: one stable crossing limit cycle with large amplitude and two unstable sliding limit cycles around each focus. See Figure \ref{fig:linear-stability}(f) for a visualization of the phase portrait with all three limit cycles. In (b) the points at which the crossing limit cycle crosses $\Sigma$ are shown as brown curves, and the maximum and minimum of the intersection of each sliding limit cycle with the (repelling) sliding region are dashed orange and purple curves, respectively. Although dynamics on $\Sigma$ are not defined \textit{a priori}, if we impose a sliding flow using the Filippov convex method, with 
\[
\lambda(y) = \frac{-a+b(y-\eta)}{-a+b(y-\eta)-p-q(y+\eta)},
\]
then there are at most two pseudoequilibria $\x_s=(x_s,y_s)$, with $x_s=0$ and $y_s$ given by the root(s) of the quadratic $\left[\lambda(y_s)(r+y_s+\eta) + (1-\lambda(y_s))(-c+y_s-\eta)\right]$. For $a=p,$ $b=q,$ and $c=r$ there is one pseudoequilibrium, at the origin.

As we decrease $\eta$ from 0.8, the two sliding cycles collide (though not at a point of tangency). Due to the repelling sliding region and nonuniqueness of solution curves leaving the sliding region, this collision leads to higher period behavior for both limit cycles while preserving the period-1 behavior, in that on the limit cycle a solution may either cycle around one focus, or both, or switch between the two. Although similar to a period-doubling bifurcation, in this case the limit cycle may have any number of periods, and the stability is unchanged by the bifurcation. This bifurcation is labeled as $(IV)$ in Figure \ref{fig:linear-stability}(b); for brevity, we refer to it as a `sliding cycle-nonunique period' bifurcation. We indicate $\eta-$values with potential higher-period limit cycles using gray shading. In reverse time, in which the unstable sliding cycles become stable and the repelling sliding region becomes attracting, this behaves like a period-doubling bifurcation of two sliding limit cycles. Neither the repelling nor attracting sliding cases appear to fit within recent classifications of piecewise-smooth bifurcations; see \cite{Kuznetsov03,di2010discontinuity,glendinning2019introduction}.

At $(III)$ in Figure \ref{fig:linear-stability}(b), the point of tangency for the right sliding limit cycle collides with the maximum of the left sliding limit cycle on $\Sigma,$ annihilating the left sliding cycle and thereby the higher-period limit cycles as well. However, the trajectory that formed the left sliding limit cycle persists, though it can no longer slide. Then at $(II)$ a similar bifurcation occurs when the minimum of the right sliding cycle collides with the tangent point of the left formerly-sliding cycle trajectory, breaking the right limit cycle. 

At $(I)$ in Figure \ref{fig:linear-stability}(b), the crossing limit cycle collides with both the left and right formerly-sliding cycle trajectories, breaking the cycle. Although the collision does not occur at the tangency points of the formerly-sliding trajectories, for brevity we refer to this as a `crossing cycle-tangency' bifurcation. As with $(IV),$ bifurcations $(I)-(III)$ do not appear to have been previously classified.

Finally, at $\eta=0$ the two points of tangency at the switching manifold collide and slide past each other, causing the repelling sliding region to become an attracting sliding region. In the context of Kuznetsov, Rinaldi, and Gragnani's classification of bifurcations in planar piecewise-smooth systems, this appears to be a $VV_2-$type bifurcation \cite{Kuznetsov03}. 

\subsection{Most probable paths}

We determine the most probable path through the switching manifold by smoothing out the boundary, then taking the limit of the resulting most probable path as our smoothing parameter goes to zero to return to the original piecewise-smooth system. 
We define the mollification of the vector field as in Equation \eqref{eq:Feps}, so that $\fscaled(x,y)=\left( f^{\e}(x,y),\, g^{\e}(x,y)  \right)=\left( f^{\e}(\x),\, g^{\e}(\x)  \right)$, where
\begin{equation}\label{eq:mollified-sys}
\begin{aligned}
f^{\e}(\x)&= \begin{cases}
\displaystyle f^-(\x) , &x\leq-\e, \\
\displaystyle f^+(\x)\Lambda\left(\frac{x}{\e}\right) +f^-(\x)\left(1-\Lambda\left(\frac{x}{\e}\right)\right) - a\Phi(-\e,x)  - p \Phi(x,\e) , &|x|<\e, \\
f^+(\x) , &x\geq\e,
\end{cases} \\
g^{\e}(\x) &= \begin{cases}
\displaystyle g^-(\x) , &x\leq-\e, \\
\displaystyle g^+(\x)\Lambda\left(\frac{x}{\e}\right) +g^-(\x)\left(1-\Lambda\left(\frac{x}{\e}\right)\right) - c\Phi(-\e,x) - r\Phi(x,\e), &|x|<\e, \\
g^+(\x) , &x\geq\e. \\
\end{cases}
\end{aligned}
\end{equation}
Here, $\Lambda(x/\e)$ is defined as in Equation \eqref{eq:Lambda}, and we define $\Phi(x_1,x_2)$ as 
\[
\Phi(x_1,x_2)=\int_{x_1}^{x_2} u \zeta_{\e}(u)\, du.
\]
Note that for $|x|<\e$, if $a=p$ the second two terms of $f^{\e}$ vanish, and if $c=r$ the second two terms of $g^{\e}$ vanish. There are now two switching manifolds in the system, $\Sigma^{\e}_{\pm} = \{\x\in\R^2:x=\pm\e\},$ but the flow is smooth across them. Note that due to the linearity of the system, mollification does not affect the fixed points $\x_{\pm}$ as long as $\e<1$. 

The persistence of the pseudoequilibrium $\x_s$ through mollification is not as straightforward as that of $\x_{\pm}$. Although a regular equilibrium may limit to a pseudoequilibrium as a smooth system becomes piecewise-smooth, nonlinear dynamics on the switching manifold may lead to multiple possible phenomena in smoothing out a piecewise-smooth field \cite{Jeffrey14}.
However, for the particular choice of $\mathbf{F}$ in this case study, after mollification there is a saddle equilibrium at the origin, which was the pre-mollification location of the pseudosaddle.

We introduce the scaling $z= x/\e,$ $z\in [-1,1]$, to simplify the limit as $\e\rightarrow 0$. With this change of variables, System \eqref{eq:Intro:DS},\eqref{eq:mollified-sys} becomes a fast-slow system,
\[
\begin{aligned}
\e \dot{z} &= f^{\e}(\e z,y), \\
\dot{y} &= g^{\e}(\e z, y).
\end{aligned}
\]
Minimizers of the rate functional $\Isig$ (Equation \eqref{eq:FW-functional}) in this rescaled flow are solutions of the EL equations
\[
\begin{aligned}
\e\ddot{z} &= \dot{y}\left( f^{\e}_y - g^{\e}_z \right) + f^{\e}f^{\e}_z + g^{\e}g^{\e}_z, \\
\ddot{y} &= \e\dot{z}\left( g^{\e}_z - f^{\e}_y \right) + f^{\e}f^{\e}_y+g^{\e}g^{\e}_y,
\end{aligned}
\]
or in Hamiltonian form, with the conjugate momenta $\varphi= \e\dot{z}-f^{\e}$ and $\psi= \dot{y}-g^{\e},$
\begin{equation}\label{eq:mollified-rescaled-sys}
\begin{aligned}
\e\dot{z} &= f^{\e}(\e z,y) + \varphi, \\
\dot{y} &= g^{\e}(\e z,y) + \psi, \\
\dot{\varphi}&=-f_z^{\e}(\e z,y)\varphi - g_z^{\e}(\e z,y)\psi,  \\
\dot{\psi}&= -f_y^{\e}(\e z,y)\varphi - \psi. \\
\end{aligned}
\end{equation}
This ensures that the stable fixed points $\x_{\pm}$ persist through the rescaling, extended by two additional components due to $\varphi$ and $\psi$. Fixed points in this extended system are given by $\textbf{z}_{\pm} = (z_{\pm},y_{\pm},0,0) = (\pm1/\e,\pm\eta,0,0)$.

\begin{figure}[t]
	\centering
	\begin{subfigure}{.34\textwidth}
		\centering
		\includegraphics[scale=1]{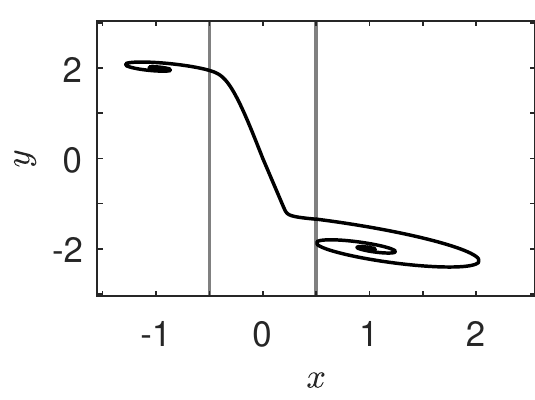}
		\caption{}
	\end{subfigure}%
	\begin{subfigure}{.33\textwidth}
		\centering
		\includegraphics[scale=1,trim=0.5cm 0cm 0cm 0cm,clip]{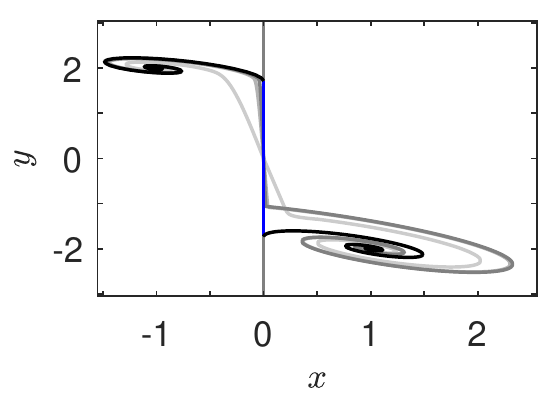}
		\caption{}
	\end{subfigure}%
	\begin{subfigure}{.33\textwidth}
		\centering
		\includegraphics[scale=1,trim=0.5cm 0cm 0cm 0cm,clip]{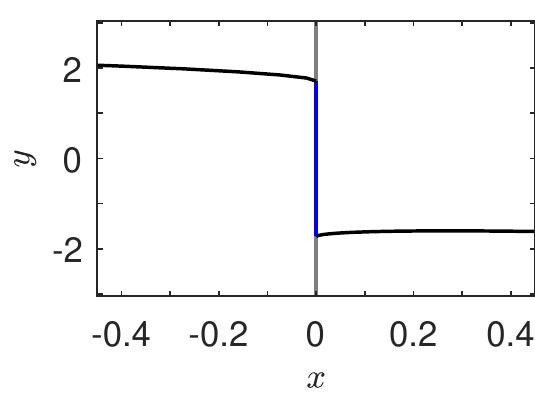}
		\caption{}
	\end{subfigure}
	\caption{(a) Most probable path of System \eqref{eq:Intro:DS},\eqref{eq:linear-system} from $\x_+$ to $\x_-$. (b) Most probable path as in (a), with varying $\e$. (c) Zoomed region of (b). All parameters given by $a=p=-2,\, b=q=-7,\, c=r=1,$ and $\eta=-2$. In (a) $\e=0.5$, and in (b),(c) $\e=0.5,\,0.1,\,0.05,$ and $\e\rightarrow 0$ (light grey to black curves).}
	\label{fig:opt-path-2}
\end{figure}

The projection onto the $xy-$plane of the solution to System \eqref{eq:mollified-rescaled-sys}, which is the most probable path of System \eqref{eq:Intro:DS},\eqref{eq:linear-system}, is shown in Figure \ref{fig:opt-path-2}. Note that we have used $x$ instead of $z$ to plot the most probable path for several values of $\e$ together. Figure \ref{fig:opt-path-2}(a) shows the most probable path for the mollified rescaled system, System \eqref{eq:mollified-rescaled-sys} with  $\e=0.5,$ including the deterministic flow from $\x_s$ to $\x_-$. Figure \ref{fig:opt-path-2}(b) shows the most probable path for decreasing values of $\e$, and Figure \ref{fig:opt-path-2}(c) shows a zoomed-in region of (b). As illustrated in (b) and (c), this is consistent with the minimizer of the rate functional for the piecewise-smooth system as $\e\rightarrow 0$. Here and in the remaining figures for this case study, we have used the function $h(x)=\exp\left(1/\left(x^2 - 1\right)\right)$ to construct the mollifier. See Appendix B for the derivation of the most probable path as the solution of System \eqref{eq:mollified-rescaled-sys}.

We can numerically validate the most probable path using Monte Carlo simulations of System \eqref{eq:Intro:SDE},\eqref{eq:linear-system} with additive noise
using the Euler-Maruyama scheme. Note that this system has a discontinuous and non-monotonic drift coefficient. Euler-Maruyama has been shown to converge in probability for such systems \cite{Gyongy96, Gyongy98}, but results for stronger convergence remain elusive; see \cite{Neuenkirch19} and references therein. Since we are concerned with convergence of the mean of the distribution of paths, Euler-Maruyama is sufficient for this investigation. 

\begin{figure}[t!]
	\centering
	\begin{subfigure}{\textwidth}
		\centering
		\includegraphics[scale=1]{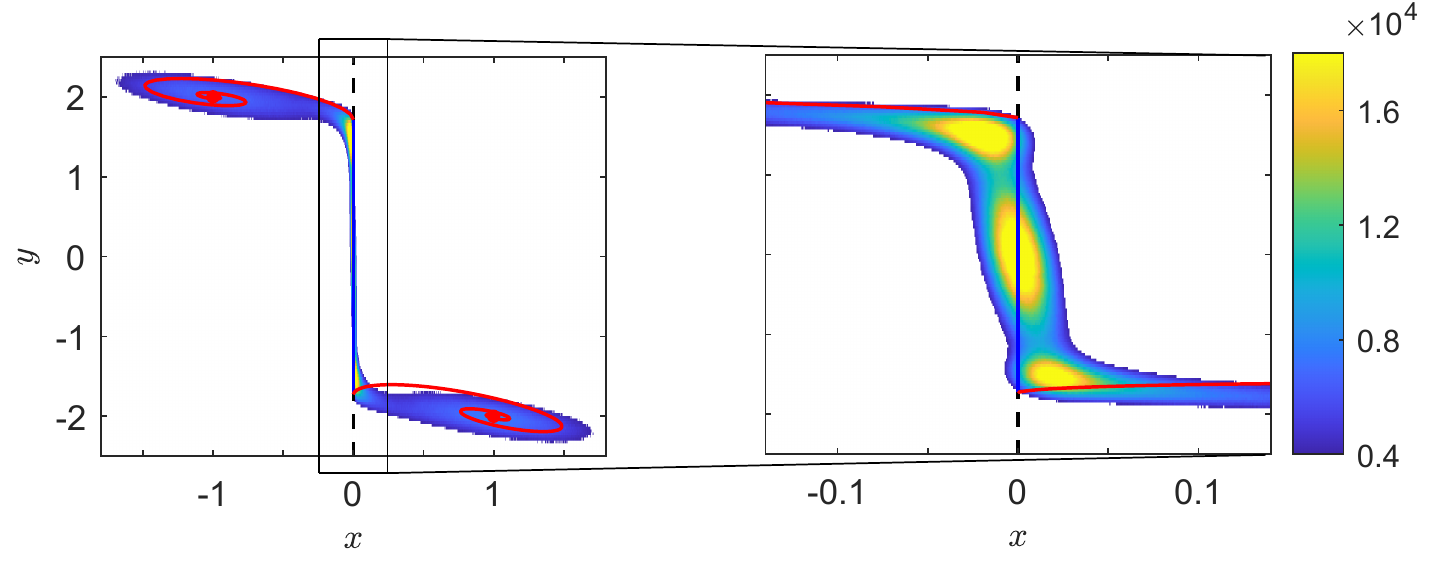}
		\caption{}
	\end{subfigure}\\
	\begin{subfigure}{\textwidth}
		\centering
		\includegraphics[scale=1]{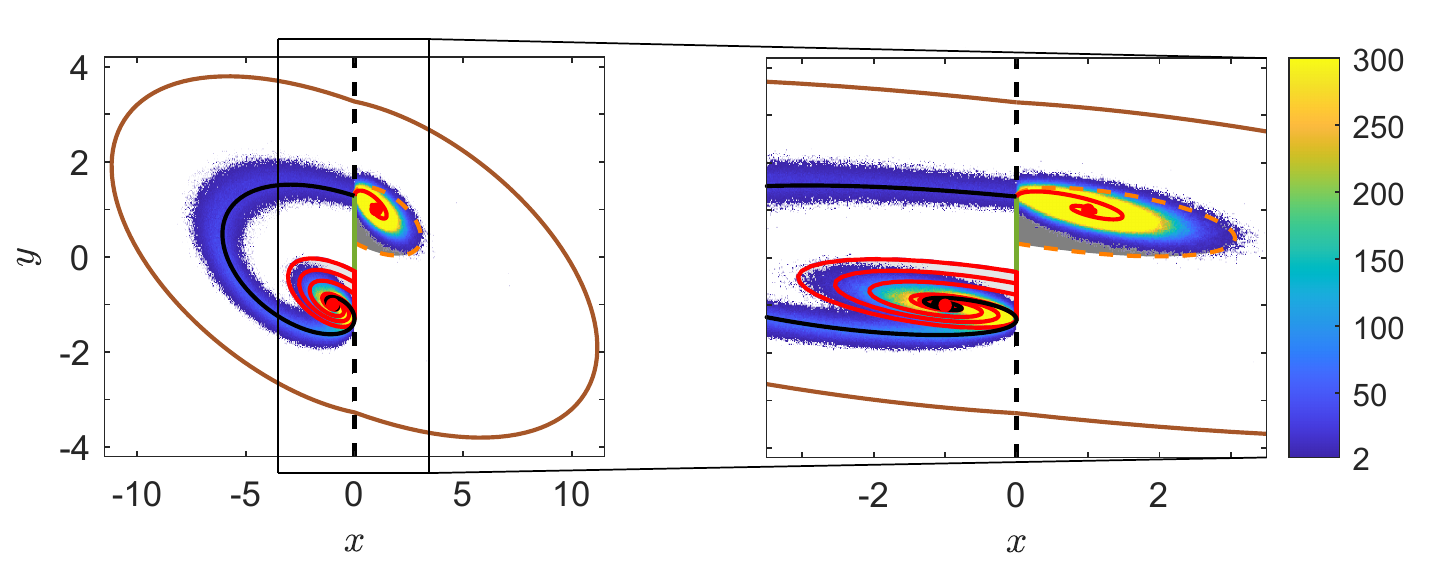}
		\caption{}
	\end{subfigure}
	
	\caption{Distribution of tipping events from Monte Carlo simulations and most probable paths (red and black curves) from $\x_+$ to $\x_-$, with (a) attracting sliding ($\eta=-2$; 6,487 tips with $N=1.056\times 10^7$ simulations) and (b) repelling sliding  ($\eta=1$; 1,960 tips with $N=1.676\times 10^7$ simulations). Parameters used are $a=p=-2,$ $b=q=-7,$ $c=r=1,$ and $\sigma=0.3$. All other curves are as in Figure \ref{fig:linear-stability}.
	}
	\label{fig:linear-Monte-Carlo}
\end{figure}

We performed Monte Carlo simulations for two values of $\eta$: $\eta = 1,$ where there is repelling sliding along $\Sigma$, and $\eta = -2,$ where there is attracting sliding. See Figure \ref{fig:linear-Monte-Carlo} for histograms of the solutions to System  \eqref{eq:Intro:SDE},\eqref{eq:linear-system} that `tip' from $\x_+$ to a neighborhood of $\x_-$. For (a) $\eta = -2$, transition paths tip over the basin boundary $\Sigma$ by tracking the attracting sliding region until it ends. Recall that in this parameter regime, there are only two limit points in the deterministic system. For (b) $\eta = 1,$ transition paths tip to a neighborhood of $\x_-$ in a sequence of basin-crossings, first over the $x>0$ unstable sliding cycle, then primarily following the deterministic flow toward the crossing limit cycle, but ultimately crossing over the $x<0$ unstable sliding cycle and approaching $\x_-$. For all figures, we have overlaid the most probable path calculated as the path $\balpha$ that minimizes the rate functional, Equation \eqref{eq:gamma-limit}. 

Note that in Figure \ref{fig:linear-Monte-Carlo}(b), we indicate both the most probable path followed by the Monte Carlo simulations (black curve) and the predicted family of nonunique most probable paths that slide (red curves). These paths coincide for $x>0$ but are distinct for $x\leq 0.$ The observed most probable path in the left-half plane corresponds to the solutions to the EL equations with boundary conditions $\x(t_1)=(0,\max_{y\in\Sigma_R} y)$ and $\x(t_f)=\x_-.$ We calculate such solutions using the gradient flow; that is, we consider solutions $\balpha:[t_1,t_f]\mapsto \R^2$ to
\begin{equation}\label{eq:grad-flow-linear}
\begin{aligned}
\frac{\partial \balpha}{\partial s} &= -\frac{\delta \Ilim}{\delta \balpha} = \ddot{\balpha} -\nabla \mathbf{F}^{\intercal}\mathbf{F} +(\nabla \mathbf{F} - \nabla  \mathbf{F}^{\intercal})\dot{\balpha},\\
\balpha(s,t_1) &= (0,\max_{y\in\Sigma_R} y), \\
\balpha(s,t_f) &= \x_-, \\
\balpha(0,t) &= \x_g(t),
\end{aligned}
\end{equation}
where $\x_g(t)$ is a smooth function satisfying $\x_g(t_1)=(0,\max_{y\in\Sigma_R} y)$ and $\x_g(t_f)=\x_-,$ and $s>0$ is an \textit{artificial time}. Along solution curves $\alpha(s,t)$, 
\[
\frac{d}{ds}\Ilim[\balpha(s,t)]\leq 0.
\]
Therefore $\Ilim$ is a Lyapunov functional and thus $\lim_{s\rightarrow \infty}\balpha(s,t)=\balpha^*(t)$, where $\balpha^*(t)$ is a solution to the EL equations, System \eqref{eq:mollified-rescaled-sys}. That is, steady state solutions of Equation \eqref{eq:grad-flow-linear} solve the Euler-Lagrange equations. The converged solution to the IBVP in Equation \eqref{eq:grad-flow-linear}, approximated using the Forward Time Centered Space finite difference scheme with $s$ as the time variable and $t$ as the space variable, is shown as the black curve in Figure \ref{fig:linear-Monte-Carlo}(b).

\begin{figure}[t!]
	\centering
	\includegraphics[scale=1]{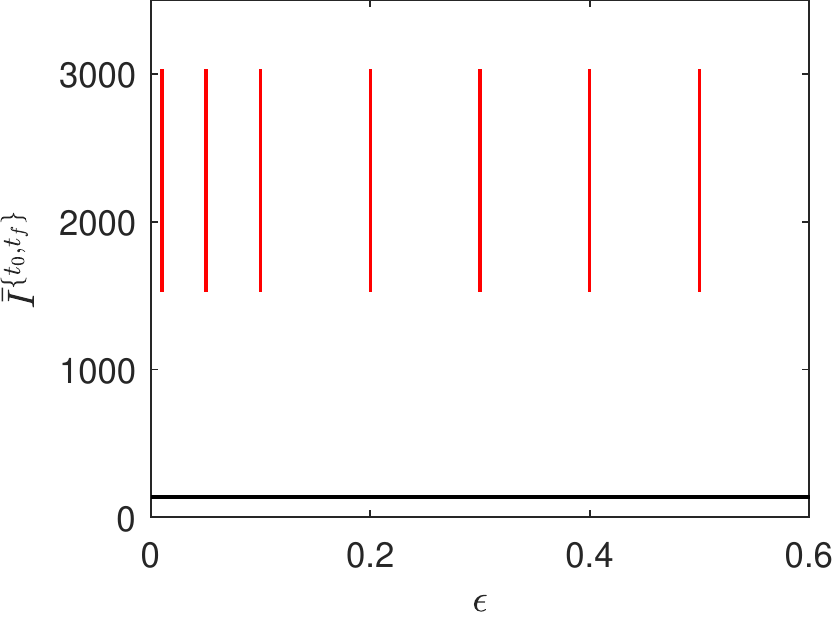}
	\caption{Values of the derived rate functional, Equation \eqref{eq:gamma-limit} for $\eta=1$ and $\mathbf{F}=\fscaled$, as $\e$ varies. The black line corresponds to the rate functional value for the crossing most probable path and the red lines correspond to the functional values for the non-unique family of sliding most probable paths.}
	\label{fig:gamma-dev}
\end{figure}

From Equation \eqref{eq:gamma-limit} we might expect the predicted most probable paths in Figure \ref{fig:linear-Monte-Carlo}(b) to minimize $\Ilim$ because $\overline{I}^{(t_1,t_f)}=0$. However, if we consider the  predicted most probable path of the piecewise-smooth system in the context of the mollified system and calculate $\IsigMol$ for several values of $\e> 0,$ the contribution of the rate functional as $\e\rightarrow 0$ is substantial: the minimum rate functional values for the predicted family of most probable paths are two magnitudes greater than those of the observed most probable path; see Figure \ref{fig:gamma-dev}. This shows anecdotally that although the minimizers of the Freidlin-Wentzell rate functional converge weakly to  minimizers of Equation \eqref{eq:gamma-limit}, they may not converge strongly. Furthermore, since the repelling sliding region has Lebesgue measure zero, the probability that the Euler-Maruyama simulations lie on $\Sigma$ at time $t$ is zero for all $t$.

\section{Case Study 2: Non-autonomous system in 1D}
\label{sec:case-study-non-autonomous}

In this section, we extend the most probable path framework to a linear one-dimensional non-autonomous system with periodic forcing. Consider Equation \eqref{eq:Intro:SDE}, where $\mathbf{F}^{\pm}:\R \times \R\mapsto \R$ are defined by 
\begin{equation}\label{eq:non-autonomous-system}
\begin{split}
\mathbf{F}^+= f^+(x,t)&=-r_{+}(x-1)+A_{+}\cos(2\pi t),\\
\mathbf{F}^-= f^-(x,t)&=-r_{-}(x-a)+A_{-}\cos(2\pi (t-p)).
\end{split}
\end{equation}
Here, $r_{\pm}, A_{\pm},p\geq 0$ and $a\leq 0$ are parameters. Since $\mathbf{F}$ is a scalar here, we will use the notation $\mathbf{F}=f$ throughout this case study. 
This case study represents a prototypical example of a one-dimensional periodically forced model with a switching manifold in the state variable, which is a typical feature of some low-dimensional conceptual climate models; see, e.g., \cite{Eisenman09,eisenman2012factors,moon2017stochastic,yang2020tipping,morupisi2021analysis,hill2016analysis}. 

Since we are considering periodic forcing, we view the flow generated by System \eqref{eq:Intro:DS},\eqref{eq:non-autonomous-system} as lying on the cylindrical extended phase space $\Pi=\R\times S^1$, where $S^1$ denotes a circle. More precisely, we make the change of variable $\tau=t$, augmenting Equation \eqref{eq:Intro:DS} by the additional equation $d \tau/dt=1$, and apply periodic boundary conditions at $\tau=0$ and $\tau=1$. However, for ease of presentation, we will continue to use Equation \eqref{eq:Intro:DS} to represent the dynamics on $\Pi$. With this convention, $\Sigma$ and $S_{\pm}$ are again defined as in Section \ref{sec:intro} and the underlying SDE is given by
\begin{equation}
dx=f(x)dt+\sigma dW, \label{Eq:NonAutSDE}
\end{equation}
where, following Filippov's convex combination method, the deterministic skeleton is given by
\begin{equation} 
\dot{x}=f(x)=\begin{cases}
f^+(x,t), & (x,t)\in S_+\cup \Sigma_+,\\
f^-(x,t), & (x,t)\in S_-\cup \Sigma_-,\\
0, & (x,t) \in \Sigma_A \cup \Sigma_R.
\end{cases}\label{Eq:DetSkeleton}
\end{equation}

\subsection{Deterministic dynamics}

Since $f(x,t)$ is linear in $x$ for $(x,t)\notin \Sigma$, the general solutions to Equation \eqref{eq:deterministic:odeFillipov} on $S_{\pm}$ are given by
\[
x_{\pm}(t)=C_{\pm}e^{-r_{\pm}t}+h_{\pm}(t),  
\]
respectively. Here $C_{\pm}\in \R$ are integration constants and
\[
\begin{split}
h_{+}(t)&=1+ \frac{r_+ A_+}{r_+^2+4\pi^2}\cos(2\pi t)+ \frac{2A_{+}\pi}{r_{+}^2+4\pi^2} \sin(2 \pi t)
,\\
h_{-}(t)&=a+ \frac{r_{-}A_{-}}{r_{-}^2+4\pi^2}\cos(2\pi (t-p))+\frac{2A_{-}\pi}{r_{-}^2+4\pi^2}\sin(2\pi(t-p)).
\end{split} 
\]
A solution with initial conditions in either $S_+$ or $S_-$ is then constructed by selecting $C_{\pm}$ to satisfy the initial condition as well as the continuity assumption on $\Sigma$. Note, in this construction solutions to Equation \eqref{eq:deterministic:odeFillipov} may be non-unique, since solutions can intersect in forward time when tracking $\Sigma_A$ and in backwards time on $\Sigma_{R}$. 

\begin{figure}[t!]
	\centering
	\begin{subfigure}[b]{0.35\textwidth}
		\includegraphics[scale=1]{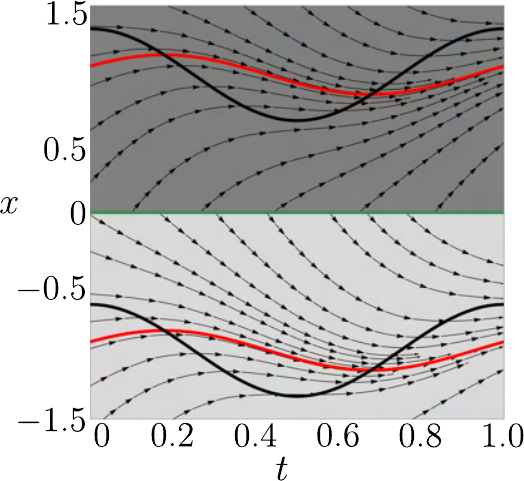}
		\caption{Neither $n_{\pm}$ intersect $\Sigma$}
	\end{subfigure}%
	\begin{subfigure}[b]{0.3\textwidth}
		\includegraphics[scale=1]{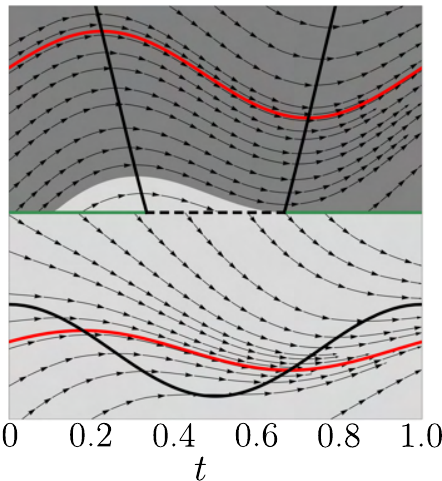}
		\caption{$n_{+}$ intersects $\Sigma$}
	\end{subfigure}%
	\begin{subfigure}[b]{0.3\textwidth}
		\includegraphics[scale=1]{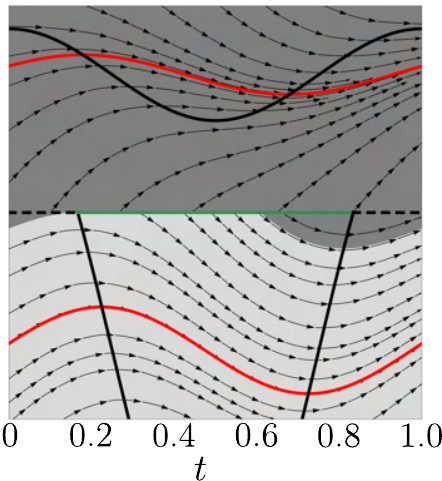}
		\caption{$n_{-}$ intersects $\Sigma$}
	\end{subfigure}\\
	~\\
	\begin{subfigure}[b]{0.35\textwidth}
		\includegraphics[scale=1]{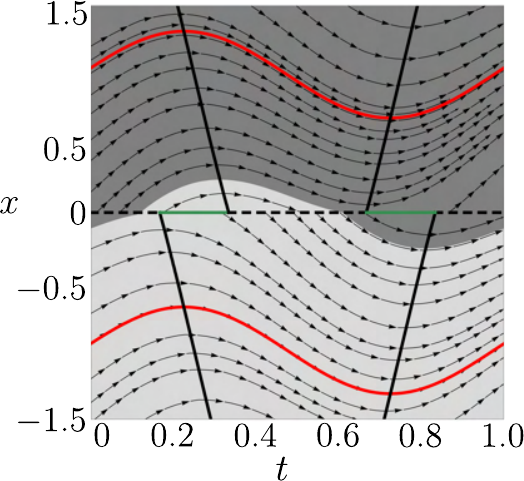}
		\caption{$n_{\pm}$ both intersect $\Sigma$}
	\end{subfigure}%
	\begin{subfigure}[b]{0.3\textwidth}
		\includegraphics[scale=1]{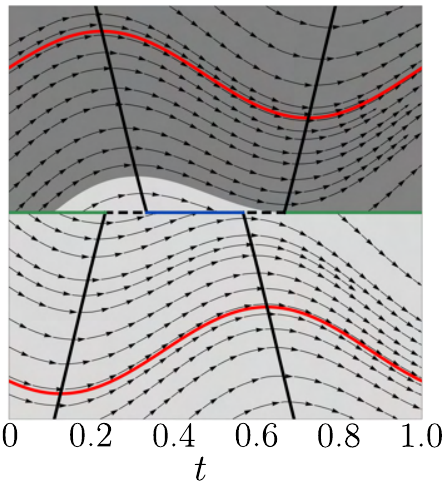}
		\caption{$n_{\pm}$ both intersect $\Sigma$}
	\end{subfigure}%
	\begin{subfigure}[b]{0.3\textwidth}
		\includegraphics[scale=1]{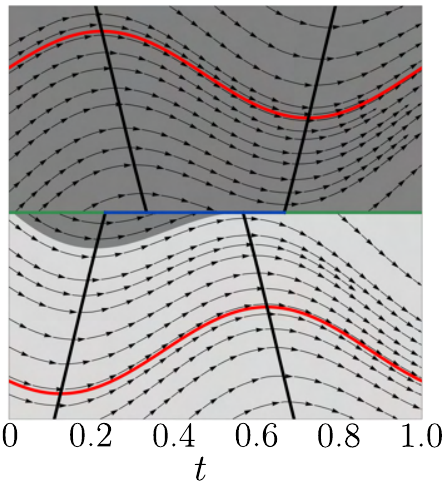}
		\caption{$n_{\pm}$ both intersect $\Sigma$}
	\end{subfigure}%
	\caption{Example phase portraits for Equations \eqref{eq:deterministic:odeFillipov}-\eqref{eq:deterministic:odeFillipovField}. The red curves show the stable limit cycles $h_{\pm}$, the black curves show the nullclines $n_{\pm}$, the green line is $\Sigma_{R}$, the blue line is $\Sigma_{A}$, and the dashed black lines are $\Sigma_{\pm}$. The dark and light gray regions correspond to the basins of attraction $B_{\pm}$ for $h_{\pm}$ respectively. 
	}
	\label{fig:Deterministic:NullclinesExample}
\end{figure}

We now consider the various parameter regimes of interest. First, if the following inequalities are satisfied:
\begin{equation}\label{Eq:Aconditions}
A_{+}< \sqrt{r_+^2+4\pi^2} \quad\text{ and }\quad A_{-}< -a \sqrt{r_{-}^2+4\pi^2}, 
\end{equation}
then $h_{\pm}(t)$ do not intersect $\Sigma$ and thus are stable limit cycle solutions to Equation \eqref{eq:deterministic:odeFillipov}. Since we are ultimately interested in noise-induced transitions from one stable limit cycle to another, we will assume these inequalities hold throughout the rest of the case study. Second, the geometry of the basins of attraction $B_{\pm}$ for $h_{\pm}(t)$, respectively, and the existence of $\Sigma_{A},$ $\Sigma_{R},$ and  $\Sigma_{\pm}$ depend on whether the nullclines of $f$ intersect $\Sigma$. The nullclines are given by curves along which $f^{\pm}(x,t)=0;$ for this system, this occurs along $n_{\pm}(t),$ where
\[
n_{+}(t)= 1+\frac{A_{+}}{r_{+}}\cos(2\pi t) \quad\mbox{ and }\quad n_{-}(t)=a+\frac{A_{-}}{r_{-}}\cos(2\pi (t-p).
\]
There are thus four cases to consider which are illustrated in Figure \ref{fig:Deterministic:NullclinesExample}:
\begin{enumerate}
	\item If $A_{+}<r_{+}$ and $A_{-}<r_{-}|a|,$ then neither nullcline intersects $\Sigma$ and thus $\Sigma_R=\Sigma$ and $\Sigma_{\pm}=\Sigma_{A}=\emptyset$; see Figure \ref{fig:Deterministic:NullclinesExample}(a). In this case $B_{\pm}$ correspond to $S_+$ and $S_-$. 
	
	\item If $A_{+}>r_{+}$ and $A_{-}<r_{-}|a|,$ then $n_+$ intersects $\Sigma$ and thus $\Sigma_- \neq \emptyset$, $\Sigma_R=\R-\Sigma_-$, and  $\Sigma_+=\Sigma_A=\emptyset$; see Figure  \ref{fig:Deterministic:NullclinesExample}(b). In this case $B_-$ has a nontrivial intersection with $S_+$ and $\Sigma$. 
	
	\item If $A_{+}<r_{+}$ and $A_{-}>r_{-}|a|,$ then $n_-$ intersects $\Sigma$ and thus $\Sigma_+\neq \emptyset$, $\Sigma_{R}=\Sigma-\Sigma_+$, and $\Sigma_-=\Sigma_A=\emptyset$; see Figure \ref{fig:Deterministic:NullclinesExample}(c). In this case $B_+$ has a nontrivial intersection with $S_-$ and $\Sigma$.
	
	\item If $A_{+}>r_{+}$ and $A_{-}>r_{-}|x_0|,$ then both $n_\pm$ intersect $\Sigma$. In this case $\Sigma_{\pm}\neq \emptyset$ but, depending on the phase shift $p$, $\Sigma_{A}$ can be either empty or nonempty; see Figure \ref{fig:Deterministic:NullclinesExample}(d)-(f). Moreover, $B_{\pm}$ can both have a nontrivial intersection with $S_+$ and $S_-$. 
\end{enumerate}

\subsection{Most probable paths}

In contrast to Case Study 1, this system does not contain stable fixed points but instead stable limit cycles. Consequently, for $t_0<t_f$, we consider the problem of determining the most probable path from $h_-(t_0)$ to $h_+(t_f)$. That is, we consider the family of optimization problems parameterized by $t_0$ and $t_f$ and define the most probable transition path from $h_-$ to $h_+$ as the minimum over this family. To do so, we redefine the admissible set of transition paths $\A$ by
\[
\A=\{\alpha \in H^{1}([t_0,t_f];\R): \alpha (t_0)=h_-(t_0) \text{ and } \alpha(t_f)=h_+(t_f)\},
\]
and consider the optimization problem
\[
\inf_{t_0<t_f}\inf_{\alpha \in \A} \Ilim[\alpha],
\]
where it follows from Theorem \ref{thm:gamma-convergence} that
\begin{equation}\label{Eq:CS2:Opt}
\Ilim[\alpha]=\int_{\I[\alpha]}\left|\dot{\alpha}(t)-f(\alpha(t),t)\right|^2\,dt+\int_{\I_{\Sigma}[\alpha]}\min_{\lambda\in [0,1]}\left\{
\left[ \lambda f^+(0,t) + \left(1-\lambda\right) f^-(0,t) \right]^2 \right\}dt.
\end{equation}

In this case study, the minimum over $\lambda\in [0,1]$ in the second integral of $\Ilim$ can be explicitly calculated.  First, for times in which $\alpha$ lies in $\Sigma_A\cup \Sigma_R$, the minimum is obtained by setting $\lambda =f^-(0,t)/\left(f^-(0,t) - f^+(0,t)\right).$ Second, for times in which $\alpha$ lies in $\Sigma_{\pm}$, $f^+(0,t)$ and $f^{-}(0,t)$  have the same sign and thus the minimum is obtained when $\lambda=0$ or $\lambda=1$. Therefore, it follows that 
\begin{equation}
\Ilim[\alpha]=\int_{\I[\alpha]}|\dot{\alpha}(t)-f(\alpha,t)|^2dt+\int_{\I_{\Sigma_{\pm}}[\alpha]}\min\{|f^+(0,t)|,|f^-(0,t)|\}^2dt. \label{Eqn:LimitFuncCase2}
\end{equation}

To study the optimization problem defined by Equation \eqref{Eqn:LimitFuncCase2} we first present the following lemma which simplifies the analysis. 
\begin{lemma}
	For the vector field defined by Equation \eqref{eq:non-autonomous-system} with parameters satisfying  \eqref{Eq:Aconditions}:
	\begin{equation}
	\inf_{t_0<t_f}\inf_{\alpha\in \mathcal{A}^{(t_0,t_f)}}\Ilim[\alpha]=\inf_{\alpha\in \mathcal{A}^{(-\infty,\infty)}}\overline{I}^{(-\infty,\infty)}[\alpha].
	\end{equation} \label{C2:InfDomain}
\end{lemma}
\begin{proof}
	Since $f^{\pm}(x,t)$ have linear growth in $x$ and are asymptotically inward-flowing, it follows from Theorem \ref{thm:min-exists} that for $t_0<t_f$ satisfying $t_0>-\infty$ and $t_f<\infty$ there exists $\alpha^*\in \mathcal{A}^{(t_0,t_f)}$ such that
	\[
	\Ilim[\alpha^*]=\inf_{\alpha \in \A}\Ilim[\alpha].
	\]
	Define $\bar{\alpha}\in \mathcal{A}^{(-\infty,\infty)}$ by
	\[
	\bar{\alpha}(t)=\begin{cases}
	h^{-}(t), & t\leq t_0,\\
	\alpha^*(t), & t_0<t<t_f,\\
	h^+(t), & t>t_f.
	\end{cases}
	\]
	By construction, $\dot{\bar{\alpha}}(t)=f(\bar{\alpha}(t),t)$ for $t\leq t_0$ and $t\geq t_f,$ and thus
	\[
	\inf_{\alpha \in \mathcal{A}^{(-\infty,\infty)}}\overline{I}^{(-\infty,\infty)}[\alpha]\leq \overline{I}^{(-\infty,\infty)}[\bar{\alpha}]=\Ilim[\alpha^*]=\inf_{\alpha \in \A}\Ilim[\alpha].
	\]
	Since this inequality is true for all values of $t_0<t_f$ satisfying $t_0\neq -\infty$ and $t_f\neq \infty,$ the result follows. 
\end{proof}

\begin{figure}[t!]
	\centering
	\begin{subfigure}[b]{0.45\textwidth}
		\centering
		\includegraphics[scale=1]{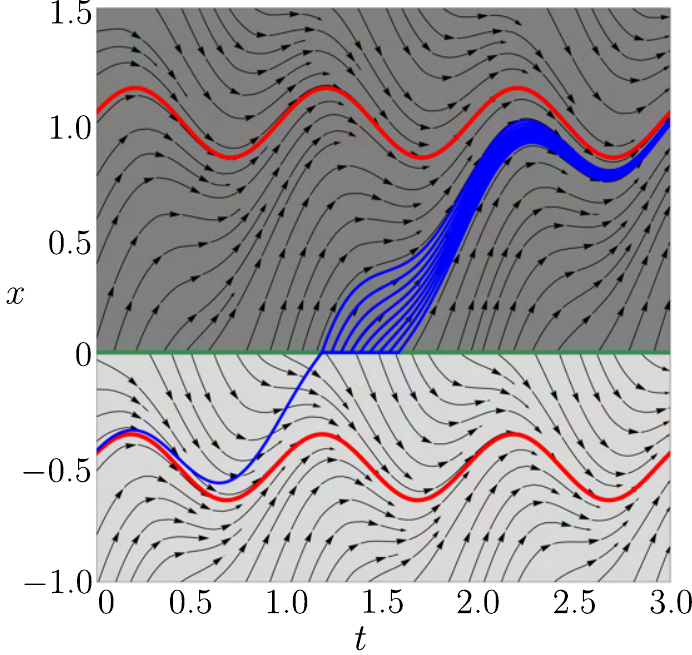}
		\caption{}
	\end{subfigure}~
	\begin{subfigure}[b]{0.5\textwidth}
		\centering
		\includegraphics[scale=1]{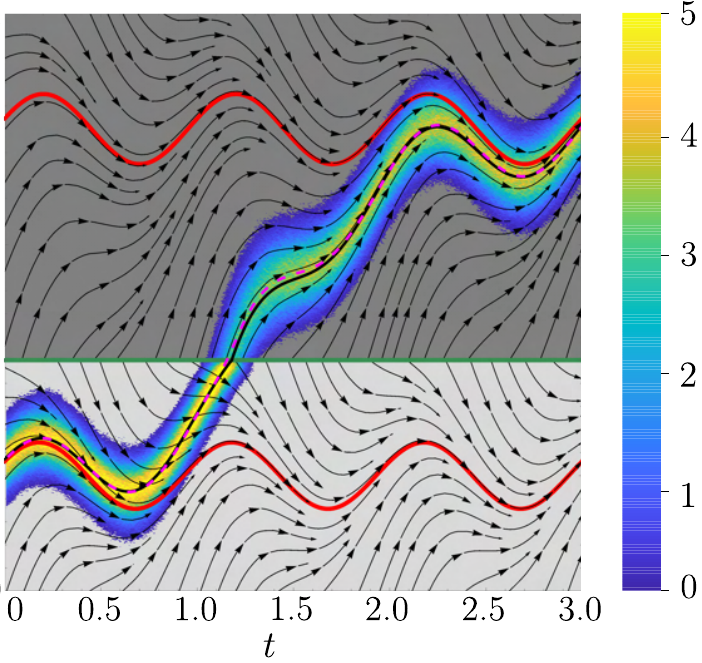}
		\caption{}
	\end{subfigure}
	\caption{(a) Non-unique minimizers of $\overline{I}^{(-\infty,\infty)}$ (blue curves) overlaid on the deterministic dynamics in the case when $\Sigma=\Sigma_R$. (b) Most probable path overlaid on the probability density generated by Monte-Carlo simulations of Equation \eqref{Eq:NonAutSDE} in a parameter regime in which $(t_{\max},0)\in B_+$. The solid black line corresponds to the most probable path constructed using Theorem \ref{Thm:NonAutMinimizer} and the dashed magenta curve is the mean of the probability density. All other curves in both (a) and (b) follow the same convention as in Figure \ref{fig:Deterministic:NullclinesExample}. The parameter values used in both (a) and (b) are $r_+ = 2,$ $r_-=3,$ $p=0$, $A_+=1,$ $A_-=1, $ and $a=-0.5$. We set $\sigma=0.2$ for the Monte-Carlo simulations.}
	\label{Fig:NonUniqueMinimizers}
\end{figure}

Interestingly, it follows from Lemma \ref{C2:InfDomain} that minimizers of Equation \eqref{Eq:CS2:Opt} are not unique. Specifically, since $f^{\pm}(x,t)=f^{\pm}(x,t+T)$ for $T\in \mathbb{Z}$, it follows that if $\alpha^*\in \mathcal{A}^{(-\infty,\infty)}$ is a minimum  then $\alpha^*_T\in \mathcal{A}^{(-\infty,\infty)}$ defined by $\alpha^*_T(t)=\alpha^*(t+T)$ is a minimum as well. Nevertheless, when represented as curves in the cylndrical phase space $\Pi$ these curves are identical. However, $\overline{I}^{(-\infty,\infty)}$ can possibly contain an uncountable number of non-unique minimizers even on $\Pi$. This possibility arises when $\Sigma_R$ and the closure of $B_+,$ which we denote as $\B_+,$ have a non-trivial intersection, and a global minimizer of $\overline{I}^{(-\infty,\infty)}$ intersects this region of $\Sigma_{R}$. Each member of the uncountable family of global minimizers can then be constructed by joining three curves as follows: the first minimizes the functional to $\Sigma_{R}$, the second tracks $\Sigma_R$ until reaching an arbitrary point $t=s$ in the intersection of $\Sigma_R$ with $\B_+$, and the third curve leaves $\Sigma_R$ and tracks the drift; see Figure \ref{Fig:NonUniqueMinimizers}(a). Key to this construction is that for Filippov systems there is a fundamental non-uniqueness in how solution curves exit a repelling sliding region, and minimizers of $\overline{I}^{(-\infty,\infty)}$ can track the drift at no cost. 

The following theorem summarizes when explicit minimizers of $\overline{I}^{(-\infty,\infty)}$ can be calculated for this case study. The key to the proof is the fact that since the system is piecewise linear, a change of variables converts the problem of determining the most probable path of the non-autonomous system to $\Sigma$ to that of an autonomous system. Moreover, the autonomous system is a gradient system and thus by Kramers' law, the most probable transition through $\Sigma$ occurs when the potential difference is a minimum \cite{berglund2013kramers, freidlin2012random}. In this case, this corresponds precisely to when the separation between $h_-$ and $\Sigma$ is minimal,
\[
t_{max}=p+\frac{1}{2\pi}\tan^{-1}\left(\frac{2\pi}{r_-}\right)+T,
\]
where $T\in\mathbb{Z};$ i.e., times when $h_-$ has a maximum.   

\begin{theorem} 
	\label{Thm:NonAutMinimizer}
	If $(t_{\max},0)\in \B_+$, where $t_{\max}=p+\frac{1}{2\pi}\tan^{-1}\left(\frac{2\pi}{r_-}\right)+T$, and $\alpha^*\in \mathcal{A}^{(-\infty,\infty)}$ is defined piecewise by
	\[
	\begin{cases}
	\alpha^*(t)=-h_-(t_{\max})e^{r_-(t-t_{\max})} +h_-(t), & t\leq t_{\max},\\
	\dot{\alpha}^*(t) \text{ satisfies \eqref{Eq:DetSkeleton}}, & t> t_{\max}, 
	\end{cases}
	\]
	then
	\[
	2h_{-}(t_{\max})^2=\overline{I}^{(-\infty,\infty)}[\alpha^*]=\inf_{\alpha\in\mathcal{A}^{(-\infty,\infty)}}\overline{I}^{(-\infty,\infty)}[\alpha].
	\]
\end{theorem}
\begin{proof}
	We first compute the most probable transition paths to $\Sigma$ by determining the most probable transition paths for the SDE
	\[
	dx=f^{-}(x,t)dt+\sigma dW
	\]
	from $h_-(t)$ to $\Sigma$. To do so, it is convenient to introduce the change of variables $y(t)=x(t)-h_-(t),$ yielding
	\[
	dy=-r_- y dt+\sigma dW,
	\]
	which corresponds to a gradient system with potential $V(y)=r_-y^2/2$. With this change of variables, the relevant admissible set $\ASigO$ and functional $\ISigO: \ASigO\mapsto \mathbb{R}$ are defined by
	\[
	\begin{aligned}
	\ASigO&=\{\gamma\in H^1([t_0,t_m];\R): \gamma(t_0)=0, \gamma(t_m)=-h_-(t_m),\text{ and } \gamma(t)<-h_-(t))\},\\
	\ISigO[\gamma]&=\int_{t_0}^{t_m}\left|\dot{\gamma}(t)
	+r_-\gamma(t)\right|^2\,dt=\int_{t_0}^{t_m}\left|\dot{\gamma}(t) + V^{\prime}(\gamma(t))\right|^2\,dt.
	\end{aligned}
	\]
	Expanding and integrating, it follows that for all $\gamma \in \ASigO$:
	\[
	\begin{aligned}
	\ISigO[\gamma]&=\int_{t_0}^{t_m}\left[\dot{\gamma}(t)^2+2\dot{\gamma}(t)V^{\prime}(\gamma(t))+V^\prime(\gamma(t))^2\right]dt \\
	&=\int_{t_0}^{t_m}\left[\left(\dot{\gamma}-V^{\prime}(\gamma(t))\right)^2+ 4\frac{d}{dt}V(\gamma(t))\right]dt\\
	&=\int_{t_0}^{t_m}\left[\dot{\gamma}-V^{\prime}(\gamma(t))\right]^2dt+4\left[V(\gamma(t_m))-V(\gamma(t_0))\right]\\
	&\geq 4V(\gamma(t_m))\\
	&=2h(t_m)^2.
	\end{aligned}
	\]
	This lower bound is an equality if $t_0=-\infty$ and $\dot{\gamma}= V^{\prime}(y)$, i.e.~$\gamma$ tracks the time-reversed dynamics. Consequently, setting  $t_m=t_{\max}$ it implies that $\gamma^*(t)=-h_-({t_{\max}})e^{r_-(t-t_{\max})}$ minimizes $\overline{I}^{(-\infty,t_m)}_{\Sigma}$ over all $t_m\in \mathbb{R}$. 
	
	Finally, for all $\alpha\in \mathcal{A}^{(-\infty,\infty)},$ if we let $t_c=\min\{t:\alpha(t)=0\}$ denote the first time $\alpha$ crosses $\Sigma$, then by the above inequality,
	\[
	\overline{I}^{(-\infty,\infty)}[\alpha]\geq \overline{I}_{\Sigma}^{(-\infty,t_c)}[\alpha-h_{-}]=2h_{-}(t_c)^2\geq 2h_{-}(t_{\max})^2.
	\]
	Therefore, since the integrand of  $\overline{I}^{(-\infty,\infty)}[\alpha]$ vanishes for times in which $\dot{\alpha}$ satisfies Equation $\eqref{Eq:DetSkeleton}$, it follows that if $(t_{\max},0)\in \B_+$ then $\alpha^*$ as constructed in the hypothesis of the theorem achieves this lower bound. 
\end{proof}

\subsection{Comparison of most probable paths to Monte-Carlo simulations}

As in Section 4, we compare the most probable paths with Monte-Carlo simulations of tipping events using the Euler-Maruyama scheme to numerically approximate realizations of Equation \eqref{Eq:NonAutSDE}. Specifically, for a realization $\chi$ of Equation \eqref{Eq:NonAutSDE} we define the following tipping time 
\[
\tau_+(\chi)=\min\{t:\chi(t)>h_+(t)\},
\]
and given $t_0<t_f$ we say $\chi$ is a tipping event on the interval $[t_0,t_f]$ if $\tau_+(\chi)<t_f$. The Monte-Carlo simulations are then conducted until the distribution of tipping events on $[t_0,t_f]$ converges. 

In the case when $(t_{\max},0)\in \B_+$, Theorem \ref{Thm:NonAutMinimizer} provides an explicit construction for the most probable paths which can be directly compared with the distribution of tipping events. In Figure \ref{Fig:NonUniqueMinimizers}(b) we plot the most probable path overlaid on the probability density of the tipping events on the interval $[0,3]$ for the same parameter values used to generate Figure \ref{Fig:NonUniqueMinimizers}(a). In this case, $\B_+=\{x:x\geq 0\}$ and thus $(t_{\max},0)\in \B_+$. Figure \ref{Fig:NonUniqueMinimizers}(b) illustrates that there is excellent agreement between the predicted most probable path and the mean and mode of the probability density generated by the Monte-Carlo simulations. Note, however, that the most probable path plotted is one that crosses and does not slide as in the other potential most probable paths illustrated in Figure \ref{Fig:NonUniqueMinimizers}(a). Clearly, while the minimizers of $\overline{I}^{(-\infty,\infty)}$ are not unique, the tipping events concentrate about the most probable path that does not slide.

When $(t_{\max},0)\notin \B_+$ the most probable paths cannot be directly computed and thus we again approximate the most probable paths by numerically computing the stationary curves of the gradient flow. In this case, to make the problem numerically tractable, we use the mollified vector field computed using a Gaussian kernel
\begin{equation}
\zeta_{\varepsilon}(x)=\frac{1}{\varepsilon \sqrt{2\pi}}\exp\left(-\frac{x^2}{2\varepsilon^2}\right).
\end{equation}
While $\zeta_{\varepsilon}$ does not have compact support, $\zeta_{1}(x)< 10^{-16}$ for $|x|\geq 9$ and thus these kernels serve as an accurate approximation of a function with compact support. With this kernel the mollified vector field for this problem becomes
\begin{equation}\label{eq:mollflow}
f^{\varepsilon}(x,t)=\varepsilon\frac{(r_{-}-r_+)}{\sqrt{2\pi}}+\frac{1}{2}\left(f_+(x,t)\left(1+\text{erf}\left(\frac{x}{\varepsilon\sqrt{2}}\right)\right)+f_{-}(x,t)\left(1-\text{erf}\left(\frac{x}{\varepsilon\sqrt{2}}\right)\right)\right),
\end{equation}
where $\text{erf}(x)=2\pi^{-\frac{1}{2}}\int_0^x\exp(-s^2)\,ds$ denotes the standard error function. The gradient flow in this case is given by the partial differential equation
\begin{equation}\label{eq:gradflow}
\begin{aligned}
\frac{\partial \alpha}{\partial s}&=-\frac{1}{2}\frac{\delta I}{\delta \alpha} =\frac{\partial^2 \alpha}{\partial t^2}-\left.\frac{\partial f^{\e}}{\partial t}\right|_{(\alpha,t)}-\left.\frac{\partial f^{\e}}{\partial \alpha}\right|_{(\alpha,t)}f^{\e}(\alpha,t),\\
\alpha(s,t_0)&=h_-(t_0), \\
\alpha(s,t_f)&=h_+(t_f), \\
\alpha(0,t)&=x_g(t),
\end{aligned}
\end{equation}
\noindent where $x_g(t)$ is a smooth function satisfying $x_g(t_0)=h_-(t_0)$ and $x_g(t_f)=h_+(t_f)$ and again $s>0$ is an \textit{artificial time}. Note, since we used a Gaussian kernel for the mollification, the functional forms of $f^{\varepsilon}$ and its various derivatives are in terms of Gaussian and error functions and thus Equation \eqref{eq:gradflow} can be numerically solved using a Forward Time Centered Space finite difference scheme. 

\begin{figure}[t!]
	\begin{subfigure}[b]{0.5\textwidth}
		\centering
		\includegraphics[scale=1]{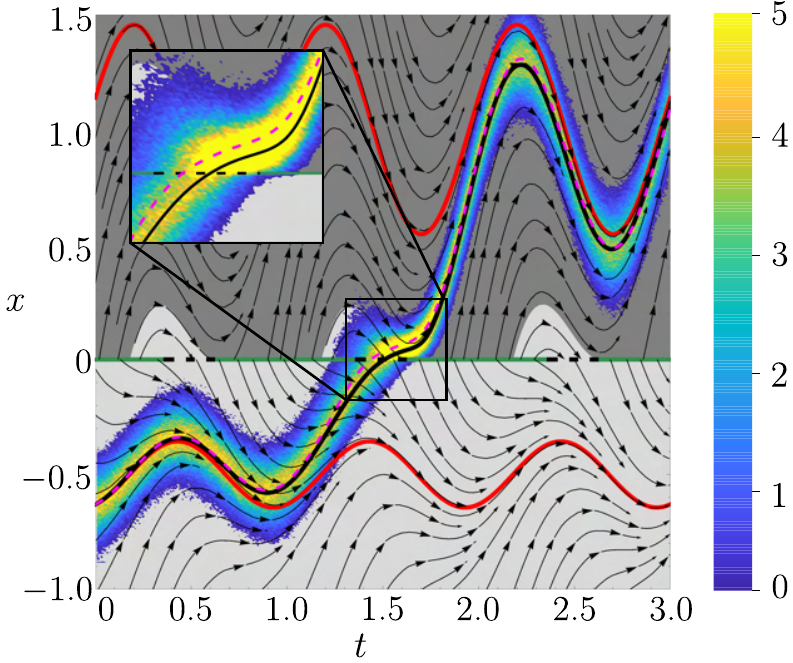}
		\caption{}
	\end{subfigure}
	~
	\begin{subfigure}[b]{0.5\textwidth}
		\centering
		\includegraphics[scale=1]{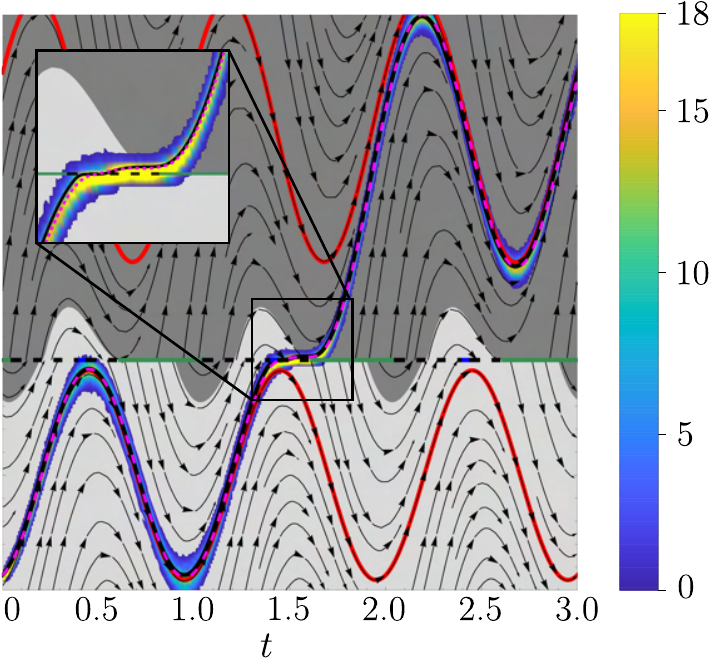}
		\caption{}
	\end{subfigure}
	\caption{Most probable paths overlaid on the probability density generated by Monte-Carlo simulations of Equation \eqref{Eq:NonAutSDE} in a parameter regime in which $(t_{\max},0)\notin B_+$. The solid black line corresponds to the most probable path numerically computed as a stationary curve of Equation \eqref{eq:gradflow} and the dashed magenta curve is the mean of the probability density. All other curves follow the same convention as in Figure \ref{Fig:NonUniqueMinimizers}. In (a) the parameters are $r_+ = 2,$ $r_-=3,$ $p=0.25,$ $A_+=3,$ $A_{-}=1,$ $a=-0.5,$ $\sigma=0.22,$ and $\e=\sigma^4$. In (b) the parameters are $r_+ = 2,$ $r_-=3,$ $p=0.25,$ $A_+=4,$ $A_{-}=3,$ $a=-0.5,$ and $\sigma=0.07,$ and $\e=\sigma^4$.}
	\label{fig:Stochastic:nullclines-2}
\end{figure}

In Figure \ref{fig:Stochastic:nullclines-2} we plot the numerically computed most probable path for the mollified vector field \eqref{eq:mollflow} overlaid on the density of tipping events for two sets of parameters in which $(t_{\max},0)\notin \B_+$. In Figure \ref{fig:Stochastic:nullclines-2}(a) the $n_+$ nullcline intersects $\Sigma$ and thus $B_-$ contains regions in which $x>0$. Moreover, $(t_{\max},0)$ lies in a crossing region. Nevertheless, in this case, the most probable path intersects $\Sigma$ slightly to the right of $(t_{\max},0)\notin \B_+$ while the mean of the density appears to pass precisely through $(t_{\max},0)$. In Figure \ref{fig:Stochastic:nullclines-2}(b) the parameters are selected so that both nullclines $n_{+}$ and $n_{-}$ intersect $\Sigma$ and thus, in addition to $B_{-}$ having nontrivial intersection with $x>0$,  $B_+$ contains regions in which $x<0$. Again, $(t_{\max},0)$ lies in a crossing region but in this case the most probable path and the mean of the tipping events are both to the left of $(t_{\max},0)$. However, in this case there is remarkable agreement between the most probable path and the mean of the tipping events. Moreover, the most probable path indicates that there is ``noise-induced'' sliding along $\Sigma$ until reaching the boundary of $B_+$.

\section{Discussion}

We have proposed an extension to the Freidlin-Wentzell theory of large deviations for determining most probable transition paths in stochastic differential equations with a piecewise-smooth drift, where the switch in drift dynamics occurs across an $(n-1)-$dimensional hyperplane. Without loss of generality, we assume this hyperplane spans the first component, so $\Sigma=\{x=0\}$. Namely, we derived an appropriate functional to determine the most probable path(s) between two metastable states on either side of $\Sigma$. For smooth regions of the drift, the derived functional matches the Friedlin-Wentzell rate functional. However, across $\Sigma$ there is an additional possible contribution to the rate functional when the most probable path slides in a crossing region of $\Sigma$.

Smoothing out the piecewise-smooth drift by mollifying it with a compactly supported, smooth, radially symmetric kernel of characteristic width $\e$ allowed us to consider most probable paths of the mollified system as minimizers of the Friedlin-Wentzell rate functional. We then considered the sequence of minimizers of the Friedlin-Wentzell rate functional in the limit as $\e\rightarrow 0.$ Using the direct method of calculus of variations, we showed that this sequence of minimizers converges weakly, and that their limit minimizes our derived functional, which itself is the $\Gamma-$limit of the sequence of rate functionals for the mollified system as $\e\rightarrow 0.$

Notably, as we recovered the piecewise-smooth drift in the limit as $\e\rightarrow 0,$ we encountered a natural choice regarding which dynamics to impose on the switching manifold $\Sigma$ and thereby determine which path across $\Sigma$ has minimal energy. Following the convention of piecewise-smooth dynamical systems, we did not choose either side, but rather introduced a scaling variable $z=x/\e\in[-1,1],$ so that as $\e\rightarrow 0$ we could still define the drift using a nonzero variable $z.$ Then the path which contributes minimal energy is tangent to $\Sigma$ in the first component and also minimizes the contribution of $z$ to the drift as $\e\rightarrow 0.$ Interpret the limiting drift on $\Sigma$ as a convex combination of the drift on either side of $\Sigma$ for a certain choice of weights, we recovered that the most probable path follows the Filippov convex combination or its time-reversed dynamics --- although we did not specify that the flow on $\Sigma$ be defined using this convention.

The derived rate functional, which has a possible contribution in addition to the Friedlin-Wentzell rate functional that depends on minimizing over some parameter $\lambda\in[0,1],$ presents an intriguing perspective on two fronts: the need for a minimum function in the rate functional and the emergence of the Filippov convex combination. The Filippov convex combination, however, is not the only possible choice for imposing a sliding flow on $\Sigma$: for example, the possibility of nonlinear sliding has recently been explored as a generalization of Filippov sliding, which may reflect observed types of sliding behavior that Filippov flow cannot reproduce; see, e.g., \cite{Jeffrey14,Novaes15}. Possible extensions of the rate functional derived here may include exploring the effect of imposing a nonlinear sliding flow in $\Sigma,$ or a complete derivation of the limiting functional for the case of a non-autonomous vector field in $\R^n.$

We used two illustrative case studies to explore implications of most probable paths in a low-dimensional setting, one case study with 2D linear piecewise-smooth drift and another with periodically-forced 1D piecewise-smooth drift. In both of these case studies, most probable paths calculated using minimizers of the derived rate functional generally matched the mean of Monte Carlo simulations of tipping events between stable states on either side of a switching manifold $\Sigma$. Additionally, both cases reveal an implication of the piecewise-smooth drift and possibility of most probable paths sliding: the most probable path is naturally non-unique when sliding in a repelling sliding region. 

However, there are still open questions arising from regimes where there is disagreement between theory and Monte Carlo simulations. Although minimizer(s) of the derived rate functional theoretically equal the most probable path only in the zero noise limit [22], we expect qualitative agreement for small noise. However, in the linear case study in Section 4, Monte Carlo simulations do not match the predicted (non-unique) family of most probable paths when the most probable paths slide in a repelling sliding region. Instead, simulations
in which tipping is observed qualitatively match the most probable path conditioned to travel from the original fixed point $\x_+$ to the onset of sliding, $(0,\max_{y\in\Sigma_R} y)$, but from the onset of sliding they follow the most probable path in the left-half plane from $(0,\max_{y\in\Sigma_R} y)$ to $\x_-$. 

Two possible explanations for this behavior are the need for a smaller noise strength, $\sigma$, or the need for a higher-order functional in the smoothing parameter, $\e$.
It may be that although the Monte Carlo simulations that tipped only happened for a small percent of simulations (0.12\% for Figure \ref{fig:linear-Monte-Carlo}(b)), the noise strength $\sigma$ needs to be smaller to capture the rare event statistics necessary to slide.

Additionally, when we consider the limiting behavior of the predicted most probable path in the piecewise-smooth limit as $\e\rightarrow 0$, although the path seen in simulations is associated with nonzero action, for $0<\e\ll 1$ the cost of sliding is at least two magnitudes greater than the cost of crossing; see Figure \ref{fig:gamma-dev} and related discussion in Section \ref{sec:case-study-2D}. 
This may indicate that in the repelling sliding case, the $\Gamma-$limit of the Friedlin-Wentzell rate functional does not completely characterize the asymptotic behavior of $\IsigMol$. In such a situation, one may consider an \textit{asymptotic development} to calculate the appropriate correction terms to the limiting functional; i.e., set
\[
\IsigMol = \overline{I}^{(t_0,t_f),0} + \e\,\overline{I}^{(t_0,t_f),1} + \e^2\, \overline{I}^{(t_0,t_f),2} + \cdots + \e^k\,\overline{I}^{(t_0,t_f),k} + o(\e^k),
\]
where $\overline{I}^{(t_0,t_f),0}=\Ilim$ is the $\Gamma-$limit of $\IsigMol$ and each $\overline{I}^{(t_0,t_f),\ell}$ for $\ell=1,2,\ldots,k$ is a functional recursively defined by the $\Gamma-$limit of $\overline{I}^{(t_0,t_f),\ell-1}$ and its infimum \cite{anzellotti1993asymptotic}. One possible extension to the present work is an asymptotic development of a family of functionals $\overline{I}^{(t_0,t_f),\ell}$, $\ell>0$, to more completely characterize the asymptotic behavior of the family $\IsigMol$.

Similarly, in the non-autonomous case study in Section \ref{sec:case-study-non-autonomous}, while the predicted most probable path may slide non-uniquely in a repelling sliding region, the mean of the observed paths does not slide. One may expect this behavior due to the sliding region having Lebesgue measure zero; see Figure \ref{Fig:NonUniqueMinimizers} and the related discussion at the end of Section \ref{sec:case-study-2D}. However, as in the linear case study, the observed most probable path intersects $\Sigma$ at the onset of the sliding region. Additionally, although the derived rate functional in Equation \eqref{eq:gamma-limit} penalizes sliding in a crossing region, in certain parameter regimes we nevertheless observe the mean of the paths exhibiting sliding behavior; see Figure \ref{fig:Stochastic:nullclines-2}(b).

The discrepancies described above for both case studies may also be related to factors not accounted for in the derived rate functional. Such factors arise in the probability density describing the transition paths $\balpha(t)$ from $\x(t_0)$ to $\x(t_f)$:
\[
P(\balpha(t))=c \exp\left[ -\frac{1}{\sigma^2} \left( \int_{t_0}^{t_f}\left\| \dot{\balpha}(t)-\mathbf{F}(\balpha(t),t)\right\|^2\,dt + \sigma^2 \int_{t_0}^{t_f} \nabla\cdot\mathbf{F}(\balpha(t))\,dt \right)\right],
\]
for some normalization factor $c\in\R$ \cite{chaichian2001path}. Note that the first integral leads to the Freidlin-Wentzell rate functional, but together both integrals are referred to as the Onsager-Machlup rate functional. Neither the normalization factor $c$ nor the Onsager-Machlup rate functional are well-defined for piecewise-smooth systems. In particular, since $\mathbf{F}$ is non-differentiable the divergence term requires considering regimes of the $\Gamma-$limit of the mollified $\fscaled$ that incorporate sigma explicitly to ensure convergence. Possible connections between most probable paths across repelling sliding regions and terms in the probability density or nonlinear convex combination will be explored in a subsequent study.

\section*{Appendix A: Inward-pointing flow and $p-$growth of $\fscaled$}
\renewcommand{\theequation}{A.\arabic{equation}}
In this appendix, we supply the proofs of Lemmas \ref{lem:p-growth} and \ref{lem:inward-Feps}.

\begin{manualtheorem}{3.1}[$p-$growth]
	There exists $\e^*$ such that for all $\e>0$ satisfying $\e<\e^*$, there exist $R^{\e}_1, c_1^{\e}>0$ and $c_2^{\e}\in \mathbb{R}$ such that $|\x|>R_1^{\e}$ implies 
	\[
	|\fscaled(\x)|\geq c_1^{\e}|\x|^p+c_2^{\e}.
	\]
\end{manualtheorem}

\begin{proof}
	For $\x=(x,\mathbf{y})$, we split the proof into the cases $|x|>\e$ and $|x|<\e$.
	
	First, assume $|x|>\e$. It follows from the reverse triangle inequality and the fact that $\zeta$ is nonnegative with unit area that
	\[
	\begin{aligned}
	\bigg| |\fscaled(\x)|-|\mathbf{F}(x)| \bigg| &\leq |\fscaled(\x)-\mathbf{F}(x)|\\
	&=\left| \int_{-\e}^{\e}\zeta_{\e}(u)\mathbf{F}(x-u,\mathbf{y})\,du-\mathbf{F}(\x)\right|\\
	&= \left| \int_{-\e}^{\e}\zeta_{\e}(u)\left(\mathbf{F}(x-u,\mathbf{y})-\mathbf{F}(\x)\right)du\right|\\
	&\leq  \int_{-\e}^{\e}\zeta_{\e}(u)\left|\mathbf{F}(x-u,\mathbf{y})-\mathbf{F}(\x)\right|du.
	\end{aligned}
	\]
	Let $R_1,c_1,c_2,c_3,c_4,c_5>0$ be defined as in Assumption \ref{ass:growth-conditions}. If $|\x|>R_1+\e$, it follows from the mean value theorem and growth condition Equation \eqref{ass:derivative-decay} that
	\[
	\begin{aligned}
	\bigg| |\fscaled(\x)|-|\mathbf{F}(x)| \bigg| &\leq \int_{-\e}^{\e} \zeta_{\e}(u)\left(c_3 \max\{|(x+\e,\mathbf{y})|^p,|(x-\e,\mathbf{y})|^p\}+c_4\right)|u|\,du\\
	&\leq\left( c_3 \max\{|(x+\e,\mathbf{y})|^p,|(x-\e,\mathbf{y})|^p\}+c_4\right)\int_{-\e}^{\e}\zeta_{\e}(u) \e du\\
	&=\e\left( c_3 \max\{|(x+\e,\mathbf{y})|^p,|(x-\e,\mathbf{y})|^p\}+c_4\right),
	\end{aligned}
	\]
	and thus
	\[
	|\fscaled(\x)|-|\mathbf{F}(\x)|>-\e\left( c_3 \max\{|(x+\e,\mathbf{y})|^p,|(x-\e,\mathbf{y})|^p\}+c_4\right).
	\]
	Consequently, since $|\x|>R_1$ and $|x|>\e$ it follows from Equation \eqref{ass:poly-growth} that 
	\[
	\begin{aligned}
	|\mathbf{F}^{\e}(\x)|&>-\e\left( c_3 \max\{|(x+\e,\mathbf{y})|^p,|(x-\e,\mathbf{y})|^p\}+c_4\right)+c_1|\x|^p+c_2\\
	&>-\e(c_3(\e+|\x|)^p+c_4)+c_1|\x|^p+c_2\\
	&>-\e\left(c_3(2|\x|)^p+c_4\right)+c_1|\x|^p+c_2\\
	&=(c_1-2^pc_3\e)|\x|^p+c_2-c_4\e.
	\end{aligned}
	\]
	
	Now, assume $|x|<\e$ and without loss of generality $x\geq 0$. Again, from the reverse triangle inequality, 
	\[
	\begin{aligned}
	\bigg| |\fscaled(\x)|-|\mathbf{F}(\x)|\bigg|\leq& \left| \int_{-\e}^{x} \zeta_{\e}(u)\left(\mathbf{F}^{+}(x-u,\mathbf{y})-\mathbf{F}^+(\x)\right)du+\int_{x}^{\e}\zeta_{\e}(u)\left(\mathbf{F}^{-}(x-u,\mathbf{y})-\mathbf{F}^{+}(\x)\right)du\right|\\
	=&\left| \int_{-\e}^{x} \zeta_{\e}(u)\left(\mathbf{F}^{+}(x-u,\mathbf{y})-\mathbf{F}^+(\x)\right)du+\int_{x}^{\e}\zeta_{\e}(u)\left(\mathbf{F}^{-}(x-u,\mathbf{y})-\mathbf{F}^{-}(0,\mathbf{y})\right)du\right.\\
	&+\left.\int_{x}^{\e}\zeta_{\e}(u)\left(\mathbf{F}^{-}(0,\mathbf{y})-\mathbf{F}^{+}(0,\mathbf{y})\right)du+\int_x^{\e}\zeta_{\e}(u)\left(\mathbf{F}^+(0,\mathbf{y})-\mathbf{F}^{+}(x,\mathbf{y})\right)du\right|\\
	\leq &  \int_{-\e}^{x} \zeta_{\e}(u)\left|\mathbf{F}^{+}(x-u,\mathbf{y})-\mathbf{F}^+(\x)\right|du+\int_{x}^{\e}\zeta_{\e}(u)\left|\mathbf{F}^{-}(x-u,\mathbf{y})-\mathbf{F}^{-}(\x)\right|du\\
	&+\left|\mathbf{F}^{-}(0,\mathbf{y})-\mathbf{F}^{+}(0,\mathbf{y})\right|+\left|\mathbf{F}^+(0,\mathbf{y})-\mathbf{F}^+(\x)\right|.
	\end{aligned}
	\]
	Therefore, if $|\x|>R_1+\e$ it follows from the mean value theorem, the growth condition Equation \eqref{ass:derivative-decay}, and the jump condition Equation \eqref{ass:jump} that
	\[
	\begin{aligned}
	\bigg| |\fscaled(\x)|-|\mathbf{F}(\x)|\bigg|\leq&  \int_{-\e}^{x} \zeta_{\e}(u)\left(c_3 |(x+\e,\mathbf{y})|^p+c_4\right)|u|du+\int_{x}^{\e}\zeta_{\e}(u)\left(c_3 |(x-\e,\mathbf{y})|^p+c_4\right)|u|du\\
	&+c_5|\mathbf{y}|^p+c_6+\e\left(c_3 |(x+\e,\mathbf{y})|^p+c_4\right)\\
	\leq& 2 \e c_4+c_3|(x+\e,\mathbf{y})|^p+c_3|(x-\e,\mathbf{y})|^p+c_5 |\mathbf{y}|^p+c_6+\e\left(c_3 |(x+\e,\mathbf{y})|^p+c_4\right).
	\end{aligned}
	\]
	Therefore,
	\[
	\begin{aligned}
	|\fscaled(\x)|&\geq  -c_5|\x|^p-c_6 -2\e c_3\max\{|2\e|^p,|\mathbf{y}|^p\}-2\e c_4+c_4+|\mathbf{F}(\x)|\\
	&\geq (c_1-c_5-2^{p+1}\e c_3)|\x|^p+(c_2-c_6-2\e c_4).
	\end{aligned}
	\]
	The same argument can be applied for $x<0$ to obtain an identical lower bound. 
	
	Therefore, from these two cases and the assumption that $c_1-c_5>0,$ if 
	\[
	\begin{aligned}
	\e&< \min\{(c_1-c_5)/(2^{p+1}c_3),c_1/(2^pc_3)\}=(c_1-c_5)/(2^{p+1}c_3),\\
	c_1^{\e}&=\min\{c_1-2^pc_3\e,c_1-c_5-2^{p+1} \e c_3\}=c_1-c_5-2^{p+1} \e c_3,\\
	c_2^{\e}&= \min\{c_2-c_4\e,c_2-c_6-2\e c_4\}=c_2-6-2\e c_4,\\
	R_1^{\e}&=R_1+\e,
	\end{aligned}
	\]
	then the result follows. 
\end{proof}

\begin{manualtheorem}{3.2}[Asymptotically inward-flowing]
	For all $\e>0$ there exist $R_2^{\e},c_7^{\e}>0$ such that $|\x|>R_2^{\e}$ implies $\fscaled(\x)\neq 0$ and
	\[
	\langle \fscaled(\x),\mathbf{r}(\x)\rangle <-c_7^{\e}|\fscaled(\x)|.
	\]
\end{manualtheorem}

\begin{proof} By construction,
	\[
	\langle \fscaled(\x),\mathbf{r}(\x)\rangle=\int_{-\e}^{\e}\zeta_{\e}(u)\langle \mathbf{F}^{-}(x-u,\mathbf{y})\mathbbm{1}_{\{u\geq x\}}(u)+\mathbf{F}^{+}(x-u,\mathbf{y})\mathbbm{1}_{\{u\leq x\}}(u),\mathbf{r}(\x)\rangle du.
	\]
	For $\x\in \mathbb{R}^n$ and $u\in (-\e, \e)$, let $\theta(u)\in [0,\pi]$ denote the angle between $\mathbf{r}(\x)$ and $\mathbf{r}(x-u,\mathbf{y})$. It follows from the mean value theorem that there exists $C(\x,u)\in (-\e, \e)$ such that
	\[
	\begin{aligned}
	|1-\cos(\theta(u))|&=\langle \mathbf{r}(\x),\mathbf{r}(\x)\rangle -\langle \mathbf{r}(\x), \mathbf{r}(x-u,\mathbf{y})\rangle \\
	&= C(\x,u)\frac{ |\mathbf{y}|^2}{|\x|\cdot|(C(\x,u)-x,\mathbf{y})|^3}u\\
	&=C(\x,u)\frac{ |\mathbf{y}|^2}{|\x|(C(\x,u)^2-2C(\x,u)x+|\x|^2)^\frac{3}{2}}u\\
	&\leq \frac{\e^2}{|\x|^2\left(1+(C(\x,u)^2-2C(\x,u)x)/|\x|^2\right)^{\frac{3}{2}}}.
	\end{aligned}
	\]
	Therefore, for all $u\in (-\e, \e)$ and $\theta^*\in (0,\pi]$ there exists $R_{\theta^*}>0$ such that $|\x|>R_{\theta^*}$ implies $\theta(\x,u)<\theta^*$. Furthermore, since $\mathbf{F}$ is asymptotically inward-flowing there exist $R_2,c_7>0$ such that if $|\x|>R_2$ then Inequality \eqref{ass:flow} is satisfied. Consequently, if $|\x|>\max\{R_2+\e,R_{\theta^*}\}$, then for $u\in(-\e, \e)$ it follows that
	\[
	\begin{aligned}
	\frac{\langle \mathbf{F}^{-}(x-u,\mathbf{y})\mathbbm{1}_{\{u\geq x\}}(u)+\mathbf{F}^{+}(x-u,\mathbf{y})\mathbbm{1}_{\{u\leq x\}}(u),\mathbf{r}(\x)\rangle }{\left|\left( \mathbf{F}^{-}(x-u,\mathbf{y})\mathbbm{1}_{\{u\geq x\}}(u)+\mathbf{F}^{+}(x-u,\mathbf{y})\mathbbm{1}_{\{u\leq x\}}(u)\right)\right|}&\leq \cos(\cos^{-1}(-c_7)+\theta^*)\\
	&=-c_7\cos(\theta^*)-\sqrt{1-c_7^2}\sin(\theta^*).
	\end{aligned}
	\]
	Hence, if we choose $\theta^*$ so that $c_7^{\e}=-c_7\cos(\theta^*)-\sqrt{1-c_7^2}\sin(\theta^*)>0$, there there exists $R_2^{\e}>0$ such that $|\x|>R_2^{\e}$ implies 
	\[
	\frac{\langle \mathbf{F}^{-}(x-u,\mathbf{y})\mathbbm{1}_{\{u\geq x\}}(u)+\mathbf{F}^{+}(x-u,\mathbf{y})\mathbbm{1}_{\{u\leq x\}}(u),\mathbf{r}(\x)\rangle }{\left|\left( \mathbf{F}^{-}(x-u,\mathbf{y})\mathbbm{1}_{\{u\geq x\}}(u)+\mathbf{F}^{+}(x-u,\mathbf{y})\mathbbm{1}_{\{u\leq x\}}(u)\right)\right|}<-c_7^{\e}.
	\]
	Therefore,
	\[
	\begin{aligned}
	\frac{\langle \fscaled(\x),\mathbf{r}(\x)\rangle}{|\fscaled(\x)|}&< -c_7^{\e}\frac{ \int_{-\e}^{\e}\zeta_{\e}(u)\left( \mathbf{F}^{-}(x-u,\mathbf{y})\mathbbm{1}_{\{u\geq x\}}(u)+\mathbf{F}^{+}(x-u,\mathbf{y})\mathbbm{1}_{\{u\leq x\}}(u)\right)du}{\left|\int_{-\e}^{\e}\zeta_{\e}(u)\left( \mathbf{F}^{-}(x-u,\mathbf{y})\mathbbm{1}_{\{u\geq x\}}(u)+\mathbf{F}^{+}(x-u,\mathbf{y})\mathbbm{1}_{\{u\leq x\}}(u)\right)du\right|}\\
	&<-c_7^{\e}\frac{ \int_{-\e}^{\e}\zeta_{\e}(u)\left( \mathbf{F}^{-}(x-u,\mathbf{y})\mathbbm{1}_{\{u\geq x\}}(u)+\mathbf{F}^{+}(x-u,\mathbf{y})\mathbbm{1}_{\{u\leq x\}}(u)\right)du}{\int_{-\e}^{\e}\zeta_{\e}(u)\left|\left( \mathbf{F}^{-}(x-u,\mathbf{y})\mathbbm{1}_{\{u\geq x\}}(u)+\mathbf{F}^{+}(x-u,\mathbf{y})\mathbbm{1}_{\{u\leq x\}}(u)\right)\right|du}\\
	&=-c_7^{\e}.
	\end{aligned}
	\]
	
\end{proof}

\section*{Appendix B: Derivation of the most probable path in Case Study~1}
\renewcommand{\theequation}{B.\arabic{equation}}

\begin{figure}[t!]
	\centering
	\includegraphics[scale=1]{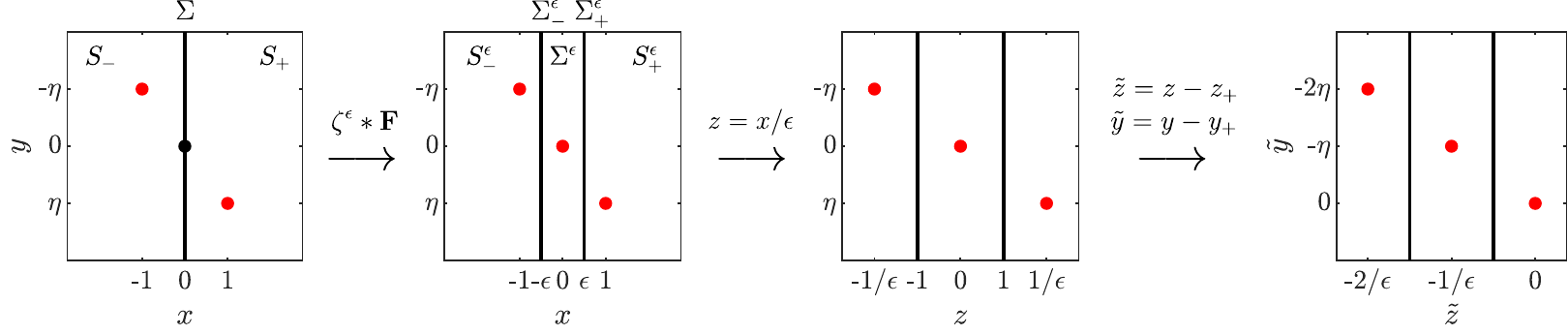}
	\caption{Diagram illustrating the transformations made to System \eqref{eq:Intro:DS},\eqref{eq:linear-system} when calculating the most probable path, for $\eta=-2$. Red points are equilibria, the black point is a pseudoequilibrium, and the vertical lines are switching manifolds.}
	\label{fig:mollification-diagram}
\end{figure}

The most probable path corresponds to solutions of System \eqref{eq:mollified-rescaled-sys} subject to appropriate boundary conditions. Since the vector field is now smooth, there must be at least one invariant manifold separating the basins of attraction of $\textbf{z}_{\pm}$. The most probable path from $\textbf{z}_+$ to $\textbf{z}_-$ follows an invariant manifold of $\textbf{z}_s$ (the pseudoequilibrium $\textbf{x}_s$ in the rescaled and extended system) in two heteroclinic orbits between $\textbf{z}_s$ and the other fixed points. We will determine the most probable path by minimizing the rate functional  $\IsigMol$ for solutions of System \eqref{eq:mollified-rescaled-sys} over two stages: \textit{(1)} from $\textbf{z}_+$ to the switching manifold $\Sigma^{\e}_+,$ and \textit{(2)} from $\Sigma^{\e}_+$ to $\textbf{z}_s$. From $\textbf{z}_s$ to $\textbf{z}_-$ the most probable path follows the deterministic flow, which does not contribute to $\IsigMol$. See Figure \ref{fig:mollification-diagram} for a diagram illustrating the sequence of mollifying, scaling, and translating transformations used in this derivation. Translating the $\textbf{z}_+$ equilibrium to the origin by defining $\tilde{z}=z-z_+$ and $\tilde{y}=y-y_+$, the boundary conditions for the path from $\mathbf{z}_+$ to $\Sigma^{\e}_+$ (Stage \textit{(1)}) become
\begin{equation}\label{eq:stage1-BCs-moll-linear}
\begin{aligned}
\tilde{z}(-\infty) &= 0, \qquad
&\tilde{y}(-\infty) &= 0, \qquad
&\varphi(-\infty) &= 0, \qquad
&\psi(-\infty) &= 0, \\
\tilde{z}(0) &= 1-\frac{1}{\e}, 
&\tilde{y}(0) &= y_1-\eta, 
&\varphi(0) &= \varphi_1,
&\psi(0) &= \psi_1.
\end{aligned}
\end{equation}
The boundary conditions $\varphi_{1}$ and $\psi_{1}$ over-determine the system if we consider that the solution must be finite at $t=-\infty$ and $t=0,$ so we allow them to vary. We will discuss choosing $y_1$ for both Stages \textit{(1)} and \textit{(2)} below. We then translate the most probable path for Stage \textit{(1)} back to the original coordinates before continuing to Stage \textit{(2)}, which gives the most probable path from $\Sigma^{\e}_+$ to $\mathbf{z}_s$. 

The boundary conditions for Stage \textit{(2)} are as follows:
\begin{equation}\label{eq:stage2-BCs-moll-linear}
\begin{aligned}
z(0) &= 1, \qquad
&y(0) &= y_1, \qquad
&\varphi(0) &= \varphi_1, \qquad
&\psi(0) &= \psi_1, \\
z(t_2) &= 0,\qquad
&y(t_2) &= 0, \qquad
&\varphi(t_2) &= 0, \qquad
&\psi(t_2) &= 0,
\end{aligned}
\end{equation}
where $t_2<\infty$ is unknown \textit{a priori}. 
Solving \eqref{eq:mollified-rescaled-sys},\eqref{eq:stage2-BCs-moll-linear} gives us the most probable path for Stage \textit{(2)}. To apply the boundary conditions, we first eliminate coefficients of the solution associated with positive eigenvalues. Applying the boundary condition at $t=0$ then leads to a system with two unknown coefficients and three initial conditions, which is over-determined. If we introduce the condition that along the most probable path the Hamiltonian of the system is zero, 
\[
H = f^{\e}\varphi + g^{\e}\psi + \frac{1}{2}\varphi^2 + \frac{1}{2}\psi^2=0,
\]
we can eliminate an equation in System \eqref{eq:mollified-rescaled-sys}. It suffices to choose the variable associated with a positive real eigenvalue and an eigenvector with nonzero conjugate momentum. For the parameters used in Figure \ref{fig:opt-path-2}, $\varphi$ satisfies this condition. There are generally two possible values of $\varphi$: we choose $\varphi= -f^{\e}+\sqrt{(f^{\e})^2-2g^{\e}\psi-\psi^2},$ so that $\varphi=0$ when $\psi=0$.

\section{Acknowledgments}

The work of JG was supported by the American Institute of Mathematics, and the work of JG and JZ was supported by Wake Forest University. Computations were performed on the Wake Forest University DEAC Cluster, a centrally managed resource with support provided in part by the University.


\bibliographystyle{unsrtnat}

\setlength{\bibsep}{2.88pt plus 0.2ex}
\urlstyle{same}

\end{document}